\documentclass[a4paper,twoside,reqno,11pt]{amsart}

\usepackage{amsmath, latexsym, amsfonts, amssymb, amsthm, amscd}
\usepackage{epsfig}
\usepackage{graphics,epsf,psfrag}

\usepackage[a4paper,scale={0.715,0.74},marginratio={1:1},footskip=7mm,headsep=10mm]{geometry}

\numberwithin{equation}{section}

\newtheorem{thm}{Theorem}

\newtheorem{prop}[thm]{Proposition}
\newtheorem{lem}[thm]{Lemma}

\newtheorem{rem}[thm]{Remark}

\def\bs{\mbox{\sf b}}

\def\l{{\ell}}

\def\NN{{\mathbb N}}

\def\PP{{\mathbb P}}

\def\G1bar{\bar{G}_1}

\def\L1bar{\bar{L}_1}

\newcommand{\bbE}{{\ensuremath{\mathbb E}} }

\newcommand{\bbP}{{\ensuremath{\mathbb P}} }

\newcommand{\cA}{{\ensuremath{\mathcal A}} }

\newcommand{\cC}{{\ensuremath{\mathcal C}} }
\newcommand{\cD}{{\ensuremath{\mathcal D}} }

\newcommand{\cJ}{{\ensuremath{\mathcal J}} }

\newcommand{\cL}{{\ensuremath{\mathcal L}} }

\newcommand{\cT}{{\ensuremath{\mathcal T}} }

\newcommand{\cZ}{{\ensuremath{\mathcal Z}} }

\newcommand{\ga}{\alpha}
\newcommand{\gb}{\beta}

\newcommand{\gd}{\delta}

\newcommand{\gs}{\sigma}

\newcommand{\go}{\omega}

\newcommand{\R}{\mathbb{R}}
\newcommand{\Z}{\mathbb{Z}}
\newcommand{\N}{\mathbb{N}}

\newcommand{\sumtwo}[2]{\sum_{\substack{#1 \\ #2}}} 

\def\bP{\ensuremath{\bs{\mathrm{P}}}}
\def\bE{\ensuremath{\bs{\mathrm{E}}}}
\def\dd{\mathrm d}
\def\cf{\textsc f}
\def\tZ{\texttt Z}
\def\tf{\texttt f}

\newcommand{\ind}{\bs{1}}

\def\bs{\boldsymbol}

\title[A diluted disordered polymer model]
{The quenched critical point\\ of a diluted disordered polymer model}

\author{Erwin Bolthausen}

\address{Institut f\"ur Mathematik, Universit\"at Z\"urich,
Winterthurerstrasse 190, CH-8057 Z\"urich}

\email{eb@math.unizh.ch}

\author{Francesco Caravenna}

\address{Dipartimento di Matematica Pura e Applicata, Universit\`a
degli Studi di Padova, via Trie\-ste 63, 35121 Padova, Italy}

\email{francesco.caravenna\@@math.unipd.it}

\author{B\'eatrice de Tili\`ere}

\address{Institut de Math\'ematiques, Universit\'e de Neuch\^atel,
Rue Emile-Argand 11, CH-2007 Neu\-ch\^a\-tel}

\email{beatrice.detiliere@unine.ch}

\date{\today}



\begin{document}

\begin{abstract}
We consider a model for a polymer interacting with an attractive wall
through a random sequence of charges. We focus on the so-called diluted
limit, when the charges are very rare but have strong intensity. In this regime,
we determine the quenched critical point of the model, showing that it
is different from the annealed one.
The proof is based on a rigorous renormalization procedure.
Applications of our results to the problem of a copolymer near a selective
interface are discussed.
\end{abstract}

\maketitle


\section{Introduction}
\label{sec:intro}

The issue addressed in this work is the determination of the
quenched critical point for the localization/delocalization
phase transition of a polymer interacting with an attractive wall
through a diluted disordered potential. The model we consider
was first introduced by Bodineau and Giacomin
in~\cite{cf:BG}, as a {\sl reduced model} for the so-called {\sl copolymer near
a selective interface model}~\cite{cf:BdH}, with the hope that it would
have the same behavior as the full copolymer model,
in the limit of weak coupling constants.
As we will see, our main result shows that this is
not the case.

The cornerstone of our approach is a
rigorous renormalization procedure. We point out that
the same result has recently been obtained by Fabio
Toninelli \cite{cf:T2}, with a rather different approach,
see the discussion following Theorem~\ref{th:main} below for more details.


\subsection{The model and the free energy}
\label{sec:model}

Let $S=\{S_n\}_{n\ge 0}$ be the simple symmetric random walk on $\Z$,
and denote by $P$ its law. More explicitly,
$S_0 = 0$ and $\{S_{n} - S_{n-1}\}_{n\ge 1}$ are i.i.d.
random variables with $P(S_1 = +1) = P(S_1 = -1) = \frac 12$.
For $N\in\N := \{1,2,\ldots\}$ we denote by
$P_N^+(\,\cdot\,)=P(\,\cdot\,|S_n\geq 0 \mbox{ for } n=1,\ldots,N)$
the law of the random walk conditioned to stay non-negative up to
time $N$. The trajectories $\{(n,S_n)\}_{0 \le n \le N}$ under $P_N^+$
model the configurations of a polymer chain of length $N$ above an
{\sl impenetrable wall}.

The interaction of the polymer with the wall
is tuned by two parameters $\beta \ge 0$, and $p \in [0,1]$.
For fixed $\gb$ and $p$,
we introduce a sequence $\omega=(\omega_n)_{n \ge 1}$ of i.i.d. random
variables, taking values in $\{0,\beta\}$ and with law $\PP$ given by:
\begin{equation} \label{eq:env}
\PP(\omega_1=\beta)=p, \quad \PP(\omega_1=0)=1-p.
\end{equation}
We are ready to define our model: for a fixed (typical) realization $\go$
and $N\in \N$, we introduce the probability measure $\bP_{N,\omega}^{\gb,p}$ defined by:
\begin{equation} \label{eq:model}
\frac{\dd\bP_{N,\omega}^{\gb,p}}{\dd P_N^+}(S)=\frac{1}{Z_{N,\omega}^{\gb,p}}
\exp\left(\sum_{n=1}^N \omega_n \ind_{\{S_n=0\}}\right),
\end{equation}
where the normalization constant $Z_{N,\go}^{\gb, p} :=
E_N^+ \big( \exp \big(\sum_{n=1}^N \omega_n \ind_{\{S_n=0\}} \big) \big)$
is usually called the {\sl partition function}.

\begin{figure}[t]
\begin{center}
\psfragscanon
\psfrag{charges}[c][c]{charges}
\psfrag{Sn}[c][c]{$S_n$}
\psfrag{n}[c][c]{$n$}
\includegraphics[width=13cm]{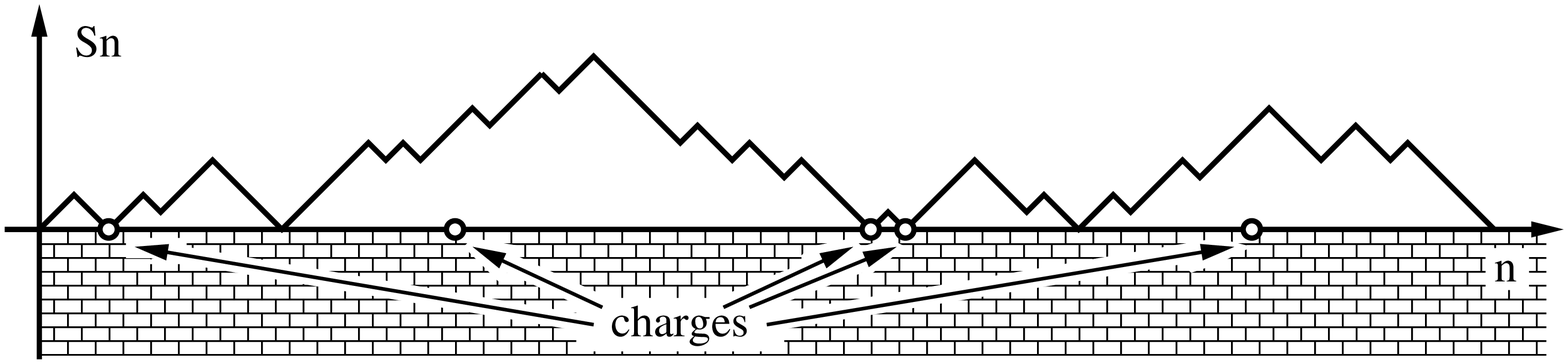}
\end{center}
\caption{A typical path of the polymer measure $\bP_{N,\go}^{\gb, c}$.}
\label{fig1}
\end{figure}


Notice that the polymer measure
$\bP_{N,\go}^{\gb, p}$ and the partition function $Z_{N,\go}^{\gb, p}$
are functions of $N$ and $\go$ only;
the superscripts $\gb, p$ are there to indicate that we are interested in the
case when the sequence $\go$ follows the law $\bbP$, which depends on $\gb, p$.

\smallskip

In this paper we focus on the regime of large $\beta$ and small $p$: then
$\go$ represents a random sequence of charges sitting on the wall (i.e.
on the $x$-axis), which are rare, but of strong intensity, and which attract
the polymer, see Figure~\ref{fig1}. We are interested in the behavior
of the polymer measure $\bP_{N,\go}^{\gb,p}$ in the limit of large~$N$:
in particular, we want to understand whether the attractive effect
of the {\sl environment} $\go$ is strong enough to pin the polymer at the
wall ({\sl localization}), or whether it is still more convenient for
the polymer to wander away from it ({\sl delocalization}), as it happens
when there are no charges.
It should be clear that we are facing a competition between energy and entropy.

\smallskip

The classical way of detecting the transition between the two
regimes mentioned above, is
to study the free energy of the model, which is defined by:
\begin{equation} \label{eq:free_energy}
f(\beta,p)=\lim_{N\rightarrow\infty}\frac{1}{N}\log Z_{N,\omega}^{\gb,p}.
\end{equation}
The existence of this limit, both $\bbP(\dd \go)$--a.s. and in $L^1(\bbP)$,
and the fact that $f(\beta,p)$ is non-random
are proved in \cite{cf:G} via super-additivity arguments.
Notice that trivially $Z_{N,\go}^{\gb,p} \ge 1$ and hence $f(\gb,p) \ge 0$
for all $\gb, p$. Zero is in fact the contribution to the free energy
of the paths that never touch the wall: indeed,
by restricting to the set of random walk trajectories that stay strictly positive
until time~$N$, one has
\begin{equation*}
    Z_{N,\go}^{\gb, p}
    \;\ge\;
    \frac{P\big(S_i > 0 \text{ for } i = 1, \ldots, N\big)}
    {P\big(S_i \ge 0 \text{ for } i = 1, \ldots, N\big)}
    \;=\;
    \frac{\frac 12 P\big(S_i \ge 0 \text{ for } i = 1, \ldots, N-1\big)}
    {P\big(S_i \ge 0 \text{ for } i = 1, \ldots, N\big)}
    \; \xrightarrow{\; N\to\infty \;} \; \frac 12\,,
\end{equation*}
where we use the well-known fact that
$P\big(S_i \ge 0 \text{ for } i = 1, \ldots, N\big) \sim (const.)/\sqrt{N}$
as $N\to\infty$, cf. \cite[Ch.~3]{cf:Fel1}.
Based on this observation,
we partition the phase space into:
\begin{itemize}
\item the Localized region $\mathcal{L} \,:=\, \{(\beta,p):\, f(\beta,p)>0\}$
\item the Delocalized region $\mathcal{D} \,:=\, \{(\beta,p):\, f(\beta,p)=0\}$.
\end{itemize}

By a standard coupling on the environment, it is clear that
the function $p \mapsto f(\gb, p )$ is non-decreasing. Therefore for every
$\gb \ge 0$, there exists a critical value $p_c(\gb) \in [0,1]$ such that the model
is localized for $p > p_c(\gb)$ and delocalized for $p < p_c(\gb)$
(in fact for $p \le p_c(\gb)$, since the function $f(\gb, p)$ is continuous).
The main goal of this work is to study the asymptotic behavior of $p_c(\gb)$,
as $\gb \to \infty$.

\smallskip
\begin{rem}\rm
One may reasonably ask to what extent the
definition of (de)local\-ization given above in terms of the free energy
corresponds to a real
(de)localized behavior of the typical paths of $\bP_{N,\go}^{\gb,p}$.
Let us just mention that, by convexity arguments, one can prove that
when $(\gb, p) \in \cL$ the typical paths of $\bP_{N,\go}^{\gb,p}$,
for large $N$, touch the wall a positive fraction of time, while
this does not happen when $(\gb, p)$ are in the interior of $\cD$.
We do not focus on path properties in this paper: for
deeper results, we refer to~\cite{cf:G}.\qed
\end{rem}

\begin{rem}\rm
Models like $\bP_{N,\go}^{\gb, p}$ are
known in the literature as (disordered) {\sl wetting models}, the terminology referring to
the interpretation of $\{S_n\}_{n\ge 0}$ as the interface of separation between
a liquid and a gazeous phase, when the liquid is above an impenetrable wall.

More generally, $\bP_{N,\go}^{\gb,p}$ belongs to the class of the so-called
{\sl disordered pinning models}, which have received a lot of attention
in the recent probabilistic literature, cf.~\cite{cf:AS,cf:GT,cf:A,cf:T,cf:T2}
(see also \cite{cf:G} for an overview).
In our case we prefer to refer to $\bP_{N,\go}^{\gb, p}$ as
a polymer model, because of its original interpretation as a simplified
model for a copolymer near a selective interface~\cite{cf:BG}
(the link with the copolymer model is discussed below).\qed
\end{rem}
\smallskip


\subsection{The main result}

\label{sec:result}

Some bounds on $p_c(\gb)$ can be obtained quite easily, as shown
in \cite[\S{}4.1]{cf:BG}. These results are stated
in the following two lemmas, whose (easy) proofs are given in detail here,
since they provide some insight into the problem. Our main result
is then stated in Theorem~\ref{th:main} below.

\smallskip
\begin{lem} \label{lem:ub}
The following relation holds:
\begin{equation}\label{eq:ub}
    - \liminf_{\gb\to\infty} \, \frac{1}{\gb} \, \log p_c(\gb) \;\le\; 1\,.
\end{equation}
\end{lem}
\smallskip

\proof
Since the limit in \eqref{eq:free_energy} holds also in $L^1(\bbP)$,
by Jensen's inequality we get
\begin{equation} \label{eq:ub1}
    f(\gb, p) \;=\; \lim_{N\to\infty} \frac 1N \, \bbE \big( \log Z_{N,\go}^{\gb, p} \big)
    \;\le \; \lim_{N\to\infty} \frac 1N \, \log \bbE \big( Z_{N,\go}^{\gb, p} \big)\,.
\end{equation}
This is usually called the {\sl annealed bound}, and the limit in the r.h.s.
above (whose existence follows by a standard super-additivity argument) is the
{\sl annealed free energy}. It can be evaluated using the
definition \eqref{eq:model} of the model, and Fubini's Theorem:
\begin{equation*}
    \bbE \big( Z_{N,\go}^{\gb, p} \big) \;=\; E_N^+ \, \bbE \,
    \Bigg( \exp\Bigg(\sum_{n=1}^N \omega_n \ind_{\{S_n=0\}}\Bigg) \Bigg)
    \; = \; E_N^+ \,
    \Bigg( \exp\Bigg(\sum_{n=1}^N \log {\rm M}(\gb, p) \ind_{\{S_n=0\}}\Bigg) \Bigg)\,,
\end{equation*}
where ${\rm M}(\gb, p) = \bbE (e^{\go_1}) = p \, e^\gb + (1-p)$.
Therefore $\bbE \big( Z_{N,\go}^{\gb, p} \big)$ is the partition
function of the simple random walk conditioned to stay non-negative, and
given a constant reward $\log {\rm M}(\gb, p)$ each time it touches zero.
This model is exactly solvable, see \cite[\S1.3]{cf:G}, and in particular we have:
\begin{equation*}
    \lim_{N\to\infty} \frac 1N \, \log \bbE \big( Z_{N,\go}^{\gb, p} \big) \;=\; 0
    \quad \qquad \text{if and only if} \quad \qquad
    {\rm M}(\gb, p) \;\le\; 2\,.
\end{equation*}
Looking back to \eqref{eq:ub1} and recalling the definition of ${\rm M}(\gb, p)$,
we have shown that
\begin{equation*}
    p \;\le\; p^a(\gb) \;:=\; \frac{1}{e^\gb - 1} \qquad \Longrightarrow \qquad (\gb, p) \in \cD\,,
\end{equation*}
where $p^a(\gb)$ is the {\sl annealed critical point}. Therefore
equation \eqref{eq:ub} is proved.\qed

\medskip

\smallskip
\begin{lem} \label{lem:lb}
The following relation holds:
\begin{equation}\label{eq:lb}
    - \limsup_{\gb\to\infty} \, \frac{1}{\gb} \, \log p_c(\gb) \;\ge\; \frac 23 \,.
\end{equation}
\end{lem}
\smallskip

\proof
We have to bound the partition
function from below. To this aim, we compute the contribution of the set
of trajectories that touch the wall wherever there is a non-zero charge
(on even sites, because of the periodicity of the random walk).
We need some notations: we introduce the subset of paths
\begin{equation*}
    \Omega_N^{\go} \;:=\; \big\{ S:\, S_n = 0 \iff \go_n > 0,\, \forall n \le N,\, n \in 2\N \big\}\,,
\end{equation*}
and the locations $\{\xi_n\}_{n\ge 0}$ of the positive charges
sitting on even sites:
\begin{equation*}
    \xi_0 \;:=\; 0 \qquad \quad
    \xi_{n+1} \;:=\; \inf \{ k > \xi_n,\, k\in 2\N : \, \go_k > 0 \}\,, \ n\in\N\,.
\end{equation*}
We denote by $\iota_N := \max\big\{ k \ge 0: \, \xi_k \le N \big\}$ the number
of positive charges (sitting on even sites) up to time $N$.
Finally, we introduce the distribution of the first return time to zero of the
simple random walk restricted to the non-negative half-line:
\begin{equation} \label{eq:K+}
    K^+(n) \;:=\; P \big( S_i > 0 \text{ for } i= 1, \ldots, n-1,\, S_n = 0 \big)
\end{equation}
(observe that $K^+(n) = 0$ for $n$ odd) and we recall that \cite[Ch.~3]{cf:Fel1}
\begin{equation} \label{eq:KSRW}
    K^+(2n) \;\overset{n\to\infty}{\sim}\; \frac{C_K}{n^{3/2}}
    \quad \text{where} \quad C_K \;=\;\frac{1}{2\sqrt{\pi}}\,, \qquad \qquad
    \sum_{n \in \N} K^+(2n) \;=\; \frac 12\,.
\end{equation}
Then we have
\begin{align*}
    Z_{N,\go}^{\gb, p} & \;\ge\; E_N^+
    \Bigg( \exp \Bigg( \sum_{n=1}^N \omega_n \ind_{\{S_n=0\}} \Bigg) \,
    \ind_{\{S \in \Omega_N^\go\}} \Bigg) \\
    & \;=\; \frac{1}{P(S_1\ge 0, \ldots, S_N \ge 0)}\cdot e^{\gb \, \iota_N}
    \cdot \Bigg( \prod_{\ell=1}^{\iota_N} K^+ \big( \xi_{\ell} - \xi_{\ell - 1} \big) \Bigg)
    \cdot \Bigg( \sum_{n = N - \xi_{\iota_N} + 1}^\infty K^+(n) \Bigg)\,.
\end{align*}
Note that $\{(\xi_\ell - \xi_{\ell - 1})/2\}_{\ell \ge 1}$ is an i.i.d. sequence of
geometric random variables with parameter $p$. Therefore, by the strong law of
large numbers, we have:
$$
\lim_{N\rightarrow\infty}\frac{\iota_N}{N}\;=\;\frac{p}{2} \ \ \text{ and } \ \
\lim_{N\rightarrow\infty}\frac{1}{N}\sum_{i=1}^{\iota_N} \log K^+(\xi_\ell -
\xi_{\ell - 1}) \;=\;
\frac{p}{2}\,\bbE \big( \log K^+(\xi_1) \big),\ \ \bbP(\dd\go)\text{--a.s.}
$$
Hence, from the last equation we get, $\bbP(\dd\go)$--a.s.,
\begin{align*}
    \lim_{N \to \infty} \, \frac 1N \, \log Z_{N,\go}^{\gb, p} \;\ge\;
    \frac p2 \, \Big( \gb \; + \; \bbE \big( \log K^+(\xi_1) \big) \Big)\,.
\end{align*}
By \eqref{eq:KSRW} there exists a positive constant $c_1$ such that
$K^+(n) \ge c_1 / n^{3/2}$, for all $n\in 2\N$.
Using this bound and Jensen's inequality yields:
\begin{align*}
    \lim_{N \to \infty} \, \frac 1N \, \log Z_{N,\go}^{\gb, p} \;\ge\;
    \frac p2 \, \Big( \gb \; + \; \log c_1 \; - \;
    \frac 32 \bbE \big( \log \xi_1 \big) \Big) \;\ge\;
    \frac p2 \, \Big( \gb \; + \; \log c_1 \; - \;
    \frac 32 \log \bbE \big( \xi_1 \big) \Big) \,.
\end{align*}
Since $\bbE \big( \xi_1 \big) = 2 p^{-1}$,
setting $c_2 := \log c_1 - \frac 32 \log 2$, we get
\begin{align*}
    \lim_{N \to \infty} \, \frac 1N \, Z_{N,\go}^{\gb, p} \;\ge\;
    \frac p2 \, \Big( \gb \; + \; \frac 32 \,\log p \;+\; c_2 \Big)\,,
\end{align*}
so that
\begin{equation*}
    p \;\ge\; e^{- \frac 23(\gb + c_2)} \qquad \Longrightarrow \qquad (\gb, p) \in \cL\,,
\end{equation*}
and equation \eqref{eq:lb} is proved.\qed

\bigskip

We can summarize Lemmas~\ref{lem:ub} and~\ref{lem:lb} in the following
way: if we knew that
\begin{equation*}
    p_c(\gb) \;\asymp\; e^{-c_{red}\, \gb} \qquad \qquad (\gb \to +\infty)\,,
\end{equation*}
then $\frac 23 \le c_{red} \le 1$
(the subscript $red$ stands for {\sl reduced model}, see the discussion below).
The main result of this paper is that in fact $c_{red} = \frac 23$. More precisely:

\smallskip
\begin{thm}\label{th:main}
For every $c> \frac 23$ there exists $\beta_0=\beta_0(c)$ such that
\begin{equation*}
    f(\gb, e^{-c\,\gb}) \;=\; 0 \qquad \text{for all } \gb \,\ge\,\gb_0\,,
\end{equation*}
i.e. $(\gb,e^{-c\,\gb}) \in \cD$ for $\gb \ge \gb_0$. Therefore
\begin{equation*}
    - \lim_{\gb \to \infty} \, \frac 1 \gb \, \log p_c(\gb) \;=\; \frac 23 \,.
\end{equation*}
\end{thm}
\smallskip

Let us discuss some consequences of this Theorem.
We recall that the model $\bP_{N,\go}^{\gb, p}$ was first introduced in
\cite{cf:BG}, as a simplified version (`reduced model')
of the so-called {\sl copolymer near
a selective interface model}, cf.~\cite{cf:BdH} (see also~\cite{cf:G} for
a recent overview). It is known that the copolymer model undergoes
a localization/delocalization phase transition. An interesting
object is the {\sl critical line} separating the two phases, in particular
in the limit of weak coupling constants, where it becomes
a straight line with positive slope $C_{cop}$.

A lot of effort has been put in finding the exact value of $C_{cop}$.
This is motivated by the fact that $C_{cop}$ appears to be a
{\sl universal} quantity: it is independent of the law of the
environment \cite[Section~3]{cf:GTdel} and
it determines the phase transition of a continuous copolymer model,
arising as the scaling limit of the discrete one \cite[\S{}0.3]{cf:BdH}.
What is known up to now is that $\frac 23 \le C_{cop} \le 1$.
Notice that $\frac 23$ and $1$ are exactly the bounds that were previously known for
$c_{red}$, and this is not by chance: indeed the definition of the model
$\bP_{N,\go}^{\gb, p}$ is inspired by the strategy behind the
proof of $C_{cop} \ge \frac 23$, cf.~\cite{cf:BG}.


The reason for introducing a reduced model
is to have a more tractable model, which would possibly have the same
behavior as the full copolymer model in the limit of weak coupling
constants, i.e. for which possibly $c_{red} = C_{cop}$.
However, the numerical results obtained in \cite{cf:CGG} provide strong
indications for the fact that $C_{cop} > \frac 23$. If this is indeed
the case, our result shows that the reduced model does not catch
the full complexity of the copolymer model, i.e. the `missing free
energy' should come from a different strategy than
the one which is at the basis of the lower bound $C_{cop} \ge \frac 23$.

By Theorem~\ref{th:main},
our model provides also a non-trivial example of a linear chain pinning model where,
for large $\gb$, the quenched critical point
$p_c(\gb)$ is different from the annealed one
$p^{a}(\gb) = 1/(e^{\gb} - 1)$ (see the proof of Lemma~\ref{lem:ub}).

\medskip

What we actually prove in this paper is a stronger version of
Theorem~\ref{th:main}, i.e. Theorem~\ref{th:main2}, stated in the next section.
The proof relies on {\sl quenched} arguments,
based on a rigorous {\sl renormalization
procedure} (somewhat in the spirit of~\cite{cf:M}). The idea
is to remove from the environment sequence $\go$ the positive charges
that are well-spaced (that therefore give no sensible contribution to
the partition function) and to cluster together the positive charges that are very close.
This procedure produces a new environment sequence
$\go'$, which has fewer charges but of stronger intensity.
The key point is that replacing $\go$ by $\go'$ in the partition function
yields an upper bound on the free energy. Then, by iterating
this transformation, we obtain environment sequences for which
the free energy can be estimated and shown to be arbitrarily small.
A detailed description of this approach, together with the organization of the paper,
is given in Section~\ref{sec:strategy}.

We point out that Theorem~\ref{th:main}
has recently been obtained by Fabio Toninelli \cite{cf:T2} with a
simpler (though more indirect) argument,
avoiding the renormalization procedure we apply. We however believe that our direct
procedure, eliminating `bad' regions in a recursive way,
is also of value for other problems, e.g. for proving that $C_{cop} < 1$.


\subsection{Beyond the simple random walk}

Theorem~\ref{th:main} can actually be extended to a broader class of models.
Namely, let $\big(\tau = \{\tau_n\}_{n\ge 0}, \bP \big)$ be a renewal
process, i.e. $\tau_0 = 0$ and $\{\tau_n - \tau_{n-1}\}_{n\ge 1}$ under $\bP$ are
i.i.d. random variables taking positive values (including $+\infty$).
It is convenient to look at $\tau$ also as the (random) subset
$\cup_{n \ge 0} \{ \tau_n\}$  of $\N\cup\{0\}$,
so that expressions like $\{k \in \tau\}$ make sense.
We assume that $\tau$ is {\sl terminating}, i.e.
\begin{equation} \label{eq:sumK}
    \gd \;:=\; \bP(\tau_1 < \infty) \;<\; 1\,,
\end{equation}
that it is {\sl aperiodic}, i.e. $\gcd \{n \in \N:\, \bP(n \in \tau) > 0 \} = 1$,
and that for some positive constant $C_K$ we have:
\begin{equation} \label{eq:asK}
    K(n) \;:=\; \bP( \tau_1 = n) \; \sim \; \frac{C_K}{n^{3/2}} \qquad \qquad (n\to\infty)\,.
\end{equation}
We introduce $\ell_N := \max \{k \ge 0:\, \tau_k \le N \}$, that gives
the number of renewal epochs up to~$N$, and the renewal function $U(\cdot)$
associated to $\tau$, defined for $n\in\N$ by
\begin{equation} \label{eq:U}
    U(n) \;:=\; \bP( n \in \tau) \;=\; \sum_{k=0}^\infty \bP( \tau_k = n )\,.
\end{equation}
By \eqref{eq:sumK} and \eqref{eq:asK}, the asymptotic behavior
of $U(n)$ is \cite[Th.~A.4]{cf:G}
\begin{equation*}
    U(n) \;\sim\; \frac{C_K}{(1-\gd)^2} \, \frac{1}{n^{3/2}} \qquad \qquad (n\to\infty)\,,
\end{equation*}
so that in particular there exists a positive constant $\cC$ such that
\begin{equation} \label{eq:boundU}
    U(n) \;\le\; \frac{\cC}{n^{3/2}} \qquad \quad \text{for every } n \in \N\,.
\end{equation}
We stress that $U(n)$ has the same polynomial behavior as $K(n)$: this is a
consequence of equation \eqref{eq:sumK} and is a crucial fact for us.

Keeping the same environment $\go=\{\go_n\}_n$ as in \eqref{eq:env},
we define the new partition function
\begin{equation} \label{eq:Znew}
    \tZ_{N,\go}^{\gb, p} \;:=\; \bE \Bigg( \exp \Bigg(
    \sum_{n=1}^{N} \go_{n} \, \ind_{\{n\in\tau\}} \Bigg) \ind_{\{N \in \tau\}} \Bigg)
    \;=\; \bE \Bigg( \exp \Bigg(
    \sum_{k=1}^{\ell_N} \go_{\tau_k} \Bigg) \ind_{\{N \in \tau\}} \Bigg)\,,
\end{equation}
and we call $\tf(\gb, p)$ the corresponding free energy:
\begin{equation} \label{eq:fnew}
    \tf(\gb, p) \;:=\; \lim_{N\to\infty} \,\frac 1N \, \log \, \tZ_{N,\go}^{\gb, p}
    \qquad \qquad \bbP(\dd\go)\text{--a.s. and in }L^1(\bbP)\,.
\end{equation}
Then we have the following extension of Theorem~\ref{th:main}.

\smallskip
\begin{thm}\label{th:main2}
For every $c> \frac 23$ there exists $\beta_0=\beta_0(c)$ such that
\rm
\begin{equation*}
    \tf (\gb, e^{-c\,\gb}) \;=\; 0 \qquad \text{for all } \gb \,\ge\,\gb_0\,.
\end{equation*}
\end{thm}
\smallskip

Let us show that Theorem~\ref{th:main} can easily be deduced from
Theorem~\ref{th:main2}. To this purpose, we choose $(\tau, \bP)$ to
be the renewal process with inter-arrival law $K^+(\cdot)$
defined in \eqref{eq:K+}. Notice that $\tau$ is not aperiodic, but
this is a minor point (it suffices to focus on the even sites to recover aperiodicity),
and that $K^+(\cdot)$ satisfies \eqref{eq:sumK} and \eqref{eq:asK} (restricted to
even sites). Then we can write the original partition function as
\begin{equation*}
    Z_{N,\go}^{\gb, p} \;=\; E_N^+ \Bigg( \exp
    \Bigg( \sum_{n=1}^{N} \ind_{\{S_n = 0\}} \Bigg) \Bigg)
    \;=\; \frac{1}{P(S_1 \ge 0, \ldots, S_N \ge 0)} \, \bE \Bigg( \exp \Bigg(
    \sum_{k=1}^{\ell_N} \go_{\tau_k} \Bigg) \Bigg)\,.
\end{equation*}
This formula looks slightly different from \eqref{eq:Znew}. First, there is
a pre-factor, due to the fact that $K^+(\cdot)$ is defined under the {\sl restricted
law} $P\big(\,\cdot\,\ind_{\{S_1 \ge 0, \ldots, S_N \ge 0\}}\big)$ while $Z_{N,\go}^{\gb, p}$
is defined as an average with respect to the {\sl conditioned law} $P_N^+$.
We have already noted that
$P\big(S_i \ge 0 \text{ for } i = 1, \ldots, N\big) \sim (const.)/\sqrt{N}$,
therefore this pre-factor is irrelevant for the purpose of computing the
free energy. The second difference is the presence in \eqref{eq:Znew} of the indicator
function $\ind_{\{N\in\tau\}}$,
but again this {\sl boundary condition} does not change
the Laplace asymptotic behavior, as is shown
in \cite[Rem.~1.2]{cf:G}. Therefore $\tf(\gb, p)$ defined by \eqref{eq:fnew} coincides
with $f(\gb,p)$ defined by \eqref{eq:free_energy}, and Theorem~\ref{th:main}
follows from Theorem~\ref{th:main2}.


\subsection{On more general return exponents}

Another natural extension is to let the renewal
process $(\tau, \bP)$ have inter-arrival distribution $K(\cdot)$ such that
\begin{equation} \label{eq:alpha}
    K(n) \;\sim\; \frac{L(n)}{n^{1+\ga}} \qquad \qquad (n\to\infty)\,,
\end{equation}
with $\ga \ge 0$ and $L: (0,\infty) \to (0, \infty)$ 
a {\sl slowly varying function} (cf.~\cite{cf:BGT}), 
i.e. such that $L(tx)/L(x) \to 1$
as $x \to \infty$, for every $t > 0$. In this paper
we stick to the case $\ga = \frac 12$ and $L(n) \to C_K$, 
for ease of notations,
but we stress that our proof of Theorem~\ref{th:main2} goes through
for the general case with easy modifications,
provided one replaces $\frac 23$ by $\frac{1}{1+\ga}$ 
in the statement of Theorem~\ref{th:main2}.
Lemma~\ref{lem:lb} also generalizes immediately, again with $\frac 23$
replaced by $\frac{1}{1+\ga}$, hence we have
$p_c(\gb) \asymp e^{-c_{red}^\ga \, \gb}$
with $c_{red}^\ga = \frac{1}{1+\ga}$.
In particular, for large $\gb$ the quenched critical point
is different from the annealed one, for any value of $\ga > 0$.


\medskip
\section{Strategy of the proof: a renormalization procedure}
\label{sec:strategy}

In this section, we explain the strategy behind the proof of
Theorem~\ref{th:main2} in details, and describe the organization of the paper.
Before doing that, we introduce some notations, and make some preliminary transformations of the partition
function $\tZ_{N,\go}^{\gb, p}$.


\subsection{Some basic notations}
\label{sec:notation}

Let us first introduce some notations to be used throughout
the paper. Given an arbitrary sequence $\go = \{\go_n\}_{n\ge 1}$, the sequence
$\{t_n\}_{n\ge 0}=\{t_n(\go)\}_{n\ge 0}$, is defined as follows:
\begin{equation}\label{eq:t}
    t_0(\go) \;:=\; 0 \qquad \qquad t_{n}(\go) \;:=\; \min\big\{ k > t_{n-1}(\go):\, \go_k > 0 \big\}\,.
\end{equation}
In our case, $\go$ is the sequence of charges, and $\{t_n\}_{n\ge 1}$
are the locations of the positive charges. However we are not assuming
that $\go$ is a typical realization of the law \eqref{eq:env}
(in particular, the size of the positive elements of $\go$ is not necessarily $\beta$).

Notice that the $\{t_n\}_{n\ge 0}$ are finite iff $\# \{i \in \N:\, \go_i > 0\} =+ \infty$,
which is always the case for us.
The increments of $\{t_n\}_{n\ge 0}$ are denoted by $\{\Delta_n\}_{n\geq 1}
=\{\Delta_n(\go)\}_{n\ge 1}$:
\begin{equation}\label{eq:Delta}
    \Delta_n(\go) \;:=\; t_n(\go) - t_{n-1}(\go) \qquad (n\in\N)\,.
\end{equation}
Finally, we introduce the sequence $\{\eta_n\}_{n\ge 1} = \{\eta_n(\go)\}_{n\ge 1}$
giving the intensities of the positive charges, i.e.
\begin{equation}\label{eq:eta}
    \eta_n(\go) \;:=\; \go_{t_n(\go)} \qquad (n\in\N)\,.
\end{equation}

We stress that one can easily reconstruct
the sequence $\go$ in terms of $\{t,\eta\} = \{t_n, \eta_n\}_n$
(or equivalently of $\{\Delta_n, \eta_n\}_n$).
Therefore we look at $\go$, $\{t, \eta\}$ and $\{\Delta, \eta\}$
as equivalent ways of describing the sequence of charges.

We also mention an elementary fact that will be used in the next section:
the sequence $\{\go_n\}_{n\ge 1}$ is i.i.d. if and only if
the variables $\{\Delta_n, \eta_n\}_{n\ge 1}$ are all independent,
the $\{\eta_n\}_{n\ge 1}$ are i.i.d. and the $\{\Delta_n\}_{n\ge 1}$
are i.i.d. geometrically distributed.


\subsection{A preliminary transformation}
\label{sec:transformation}

In the partition function $\tZ_{N,\go}^{\gb, p}$, the parameter $N$
represents the size of the system. However it turns out to be more
convenient to consider a partition function with a {\sl fixed} number
of positive charges.

Let us be more precise.
Since the limit in \eqref{eq:fnew} exists $\bbP(\dd\go)$--a.s., we can take
it along the (random) subsequence $\{t_n\}_{n\ge 0} = \{t_n(\go)\}_{n\ge 0}$, i.e. we can write,
$\bbP(\dd\go)$--a.s.,
\begin{equation} \label{eq:f2}
    \tf(\gb, p) \;=\; \lim_{n\to\infty} \, \frac{1}{t_n} \, \log
    \tZ_{t_n, \go}^{\gb,p} \,,
\end{equation}
Let us focus on $\tZ_{t_n, \go}^{\gb,p}$: by summing over the locations
of the positive charges that are visited and recalling the definition \eqref{eq:U}
of the renewal function $U(\cdot)$, we obtain the explicit expression
\begin{equation} \label{eq:boundZ}
    \tZ_{t_n,\go}^{\gb, p} \;\le\; \sum_{k=1}^{n} \sumtwo{j_1, \ldots, j_{k-1} \in \N}
    {0 =: j_0 < j_1 < \ldots < j_{k-1} < j_k = n} \prod_{\ell=1}^k
  e^{\eta_{j_\ell}} \cdot U(t_{j_\ell} - t_{j_{\ell-1}}) \,.
\end{equation}
We stress that in the r.h.s. we have $U(t_{j_\ell} - t_{j_{\ell-1}})$ and
not $K(t_{j_\ell} - t_{j_{\ell-1}})$: in fact we are fixing the positive
charges that are visited, but the path is still free to touch
the wall in between the positive charges. Also notice that,
when $\go$ is distributed according to \eqref{eq:env}, we have
$\eta_{j_\ell} = \beta$, but we keep the notation
implicit for later convenience.

\smallskip

Formula \eqref{eq:boundZ}
leads us to the following definition: for $n \in \N$ and $C \in \R^+$ we set:
\begin{equation} \label{eq:start}
    \cZ_{n}(\go, C) \; := \; \sum_{k=1}^{n} \sumtwo{j_1, \ldots, j_{k-1} \in \N}
    {0 =: j_0 < j_1 < \ldots < j_{k-1} < j_k := n} \prod_{\ell=1}^k
    e^{\eta_{j_\ell}} \cdot \frac{C}{(t_{j_\ell} - t_{j_{\ell-1}})^{3/2}} \,,
\end{equation}
so that applying the upper bound \eqref{eq:boundU} on $U(\cdot)$ we have
\begin{equation}\label{eq:miao}
\tZ_{t_n,\go}^{\gb, p} \;\le\; \cZ_{n}(\go, \cC)\,.
\end{equation}
Notice that $\cZ_n(\go, \cC)$ carries no explicit dependence of
the parameters $\gb, p$.
In fact, we look at $\cZ_{n}(\go, \cC)$ as a deterministic function of the constant $\cC$
and of the sequence $\go$. In the sequel, $\cC$ will always denote the constant
appearing in the r.h.s. of \eqref{eq:boundU}.

Now assume that $(\go = \{\go_n\}_n, \bbP_\mu)$ is an i.i.d. sequence
of random variables
with marginal law $\mu$: $\bbP_\mu (\go_1 \in \dd x) = \mu(\dd x)$. Then
for $C>0$ we define the free energy $\cf(\mu,C)$ as the limit
\begin{equation}\label{eq:flast}
    \cf(\mu,C) \;:=\; \lim_{n\to\infty} \, \frac{1}{t_n(\go)} \, \log \cZ_{n}(\go, C)
    \qquad \quad \bbP_\mu (\dd \go)\text{--a.s.}\,.
\end{equation}
If we denote by $\mu_{\gb, p} := (1-p) \, \gd_{\{0\}} + p \, \gd_{\{\gb\}}$, then by \eqref{eq:f2}
and \eqref{eq:miao} we can write
\begin{equation*}
    \tf (\gb,p) \;\le\; \cf\big( \mu_{\gb, p}\,, \cC \big)\,.
\end{equation*}
Then to prove Theorem~\ref{th:main2} it suffices to prove the following:
\emph{for every $\cC > 0$ and $c > \frac 23$ there exists $\gb_0 = \gb_0(\cC, c)$
such that $\cf\big( \mu_{\gb, \exp(-c\,\gb)}\,, \cC \big) = 0$ for every $\gb \ge \gb_0$.}

\smallskip

Therefore, from now on, we focus on $\cZ_{n}(\go, \cC)$ and $\cf(\mu, \cC)$.
The constants $\cC > 0$ and $c > \frac 23$ are fixed
throughout the sequel. We also use the shorthand $\mu_\gb$
for $\mu_{\gb, \exp(-c\,\gb)}$, i.e.
\begin{equation} \label{eq:mubeta}
    \mu_\gb \;:=\; \big( 1-e^{-c\,\gb} \big) \, \gd_{\{0\}} + e^{-c\,\gb} \, \gd_{\{\gb\}}\,.
\end{equation}
Our goal is to show that $\cf(\mu_\gb, \cC) = 0$ for $\gb$ large.
For ease of notation, {\sl we only consider the case $\gb \in \N$}
(the general case can be recovered with minor modifications).


\subsection{The renormalization procedure}
\label{sec:renormalization}

The proof of Theorem~\ref{th:main2} is achieved through an induction
argument. The steps of the induction are labeled by $\{\gb, \gb+1, \gb+2, \ldots\}$,
and we call them {\sl level $\gb$}, {\sl level $\gb + 1$}, \ldots

Each induction step consists of a renormalization procedure
that acts both on the environment $\go$, and on the partition function
$\cZ_{n}(\go, \cC)$, and produces an upper bound on the free energy $\cf(\mu, \cC)$.
Let us be more precise, by describing in detail how this
procedure works.

\medskip

\noindent
{\it Renormalizing the environment (Section~\ref{sec:env}).}
At the starting point (level $\gb$) the environment $\go$
is i.i.d. with marginal law $\mu_\gb$ defined in \eqref{eq:mubeta},
supported on $\{0\} \cup \{\gb\}$.
More generally, at level $b \ge \gb$ the environment $\go$
will be i.i.d. with marginal law $\mu_b$ supported on $\{0\} \cup \{b, b+1, \ldots\}$.
If we are at level $b$, we define a renormalization map $T_b$
acting on $\go$ that produces a new sequence $T_b(\go)$ as follows.

\begin{figure}[h]
\begin{center}
\psfragscanon
\psfrag{good charge}[c][c]{good charge}
\psfrag{isolated charges}[c][c]{isolated charges}
\psfrag{b}[c][c]{\small $b$}
\psfrag{b+5}[c][c]{\small $b+5$}
\psfrag{b+8}[c][c]{\small $b+8$}
\psfrag{Lb=}[c][c]{$L_b=$}
\psfrag{bad blocks}[c][c]{bad blocks}
\includegraphics[width=13cm]{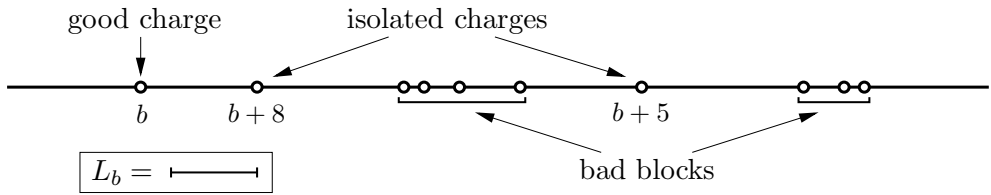}
\end{center}
\caption{{\sl Good charges}, {\sl isolated charges} and
{\sl bad blocks} at level $b$.}
\label{fig2}
\end{figure}

We first need to define {\sl isolated charges}, {\sl good charges} and {\sl bad blocks}
at level $b$. To this purpose, we fix the threshold
$L_b := \big\lfloor e^{\frac{2}{3}(b+K_b)} \big\rfloor$, where the constant
$K_b$ is defined explicitly in \eqref{eq:Kb} below. A positive charge is said to be an {\sl isolated charge}
if both its neighboring positive charges are at distance greater than $L_b$.
Among the isolated charges, we call {\sl good charges} those that have
intensity exactly equal to $b$, i.e. the least possible intensity.
Finally, a group of adjacent positive charges is said to be a {\sl bad block} if
all the distances between neighboring charges inside the group are smaller than $L_b$.
Note that a charge is either isolated, or it belongs to a bad block
(see Figure~\ref{fig2} for a graphical illustration).

Then the renormalized environment $\go' = T_b(\go)$ is obtained
from $\go$ in the following way:
each {\sl bad block} is clustered into one single larger charge,
each {\sl good charge} is erased, the {\sl isolated charges}
that are not good are left unchanged and finally
the distances between charges are suitably shortened.
In Section~\ref{sec:env} we show that the new environment $\go'$
constructed in this way is still i.i.d.
and we obtain an explicit expression for the marginal law
of $\go'_1$, denoted by $\mu_{b+1} := T_b(\mu_b)$.
Observe that by construction $\mu_{b+1}$ is supported
by $\{0\} \cup \{b+1, b+2, \ldots\}$.

\bigskip

\noindent
{\it Renormalizing the partition function (Section~\ref{sec:pf}).}
The idea behind the definition of good charges and bad blocks is the following:
\begin{itemize}
\item[-] if a charge is good, it is not worth for the polymer to visit it, because this
would entail a substantial entropy loss;
\item[-] on the other hand, if a charge belongs to a bad block
and the polymer visits it, it
is extremely convenient for the polymer to visit all the charges in the block.
\end{itemize}
These rough considerations are made precise in Section~\ref{sec:pf}
(as a side remark, notice that the choice of
what is `good' and what is `bad' is biased by the fact that we aim at proving
delocalization).
As a consequence, if we replace the environment $\go$ with the renormalized one $T_b(\go)$,
we get an upper bound on the partition function. More precisely, if we also
denote by $T_b$ the transformation acting
on $C>0$ by $T_b(C) := C \cdot (1 + B\,e^{-K_b}\, C)$
(where $B$ is an absolute constant defined in Lemma~\ref{lem:ren_th}
below), then
we show that the partition function satisfies, for every $N\in\N$,
\begin{equation} \label{eq:sloppy}
    \cZ_{n}(\go, C) \;\le\; (const.) \, \cZ_{N} \big( T_b(\go), T_b(C)\big)\,,
\end{equation}
for a suitable $n = n(\go, N)$ such that $n \ge N$ and $t_n(\go) \ge t_N(T_b(\go))$.
Taking `$\frac{1}{t_n} \log$' on both sides
of \eqref{eq:sloppy}, and letting $n\to\infty$, we then obtain for every $b \ge \gb$,
\begin{equation*}
    \cf \big(\mu_b, C) \;\le\; \cf \big( \mu_{b+1}, T_b(C) \big)
\end{equation*}
(recall that $\mu_{b+1} = T_b(\mu_b)$) and by iteration we have for $b\ge \gb$
\begin{equation} \label{eq:iter_bound}
    \cf \big(\mu_\gb, \cC) \;\le\; \cf \big( \mu_{b}, \cC_b \big) \qquad \text{where}
    \qquad \cC_b \;:=\; \big( T_{b-1} \cdot T_{b-2} \cdots T_\gb \big) (\cC)\,.
\end{equation}

\bigskip

\noindent
{\it Completing the proof (Section~\ref{sec:proof}).}
The last step is to get a control on the law $\mu_b$ and on the constant $\cC_b$,
in order to extract explicit bounds from \eqref{eq:iter_bound}.
By easy estimates, we show that $\cC_b \,\le\, 2\cC$ for every $b$,
so that this yields no problem. The crucial point is rather in estimating the law
$\mu_b$: in Proposition~\ref{propInduction} we prove
(when $\gb$ is large but fixed) an explicit stochastic domination of $\mu_b$,
which allows to show that
\begin{equation*}
    \lim_{b \to\infty} \, \cf \big(\mu_b, \cC_b) \; = \; 0 \,.
\end{equation*}
By \eqref{eq:iter_bound} this implies that
$\cf \big(\mu_\gb, \cC) = 0$, and Theorem~\ref{th:main2} is proved.


\medskip
\section{Renormalization of the environment}
\label{sec:env}

In this section we describe the renormalization transformation performed on
the environment, outlined in the
previous section. At level
$b\in\{\beta,\beta+1,\ldots\}$ the sequence $\omega=\left\{  \omega_{1}%
,\omega_{2},\ldots\right\}  $ is i.i.d. where the $\omega_{i}$ have law
$\mu_{b}$ supported on $\{0\}\cup\{b,b+1,\ldots\}$. We set%
\begin{equation}
c_{b}\;:=\; \mu_{b}\left(  [b,\infty)\right) \,,
\qquad \ \ \
\tilde c_b \;:=\; \mu_{b} \left( [b+1, \infty) \right) \,.
\label{cb-Definition}%
\end{equation}
Then we act on $\omega$ by clustering and removing charges, in order to obtain
a new environment sequence $\go' = T_{b}(\omega)$ which is still i.i.d. but with a
new charge distribution $\mu_{b+1}=T_{b}(\mu_{b})$ supported by $\{0\}\cup
\{b+1,b+2,\ldots\}$. The definition of $T_b$ depends on%
\begin{equation}
K_{b}:=\big\lfloor K_{0}+\log^+(2\mathcal{C})+2\log b\big\rfloor
\,,\qquad L_{b}:=\big\lfloor \mathrm{e}^{\frac{2}{3}(b+K_{b})}\big\rfloor\,,
\label{eq:Kb}%
\end{equation}
where $\log^+(x) := \log(x) \vee 0$, $K_{0}$ is an absolute positive constant
defined in Lemma~\ref{lem:ren_th} below, and
$\mathcal{C}$ is the constant appearing in \eqref{eq:boundU}.

We remind the reader that the sequence $\left\{\omega_{k}\right\}  _{k\geq
1}$ is in one-to-one correspondence with the pair of sequences $\big(
\left\{  t_{k}\right\}  _{k\geq 1},\left\{  \eta_{k}\right\}  _{k\geq 1}\big)$,
where $0<t_{1}<t_{2}<\cdots$ is the sequence of successive times where
$\omega_{t_{k}}>0$ and $\eta_{k}:=\omega_{t_{k}}$. We also set for
convenience $\eta_0 := 0$. Alternatively, we can
replace $\left\{  t_{k}\right\}  _{k\geq1}$ by the sequence $\left\{
\Delta_{k}\right\}  _{k\geq1},\ \Delta_{k}:=t_{k}-t_{k-1},$ where $t_{0}:=0.$
We will freely switch from one representation to the other without special mentioning
(see \S\ref{sec:notation} for more details).

\begin{figure}[t]
\begin{center}
\psfragscanon
\psfrag{intensity b = good charge}[c][c]{intensity $b$ $=$ good charge}
\psfrag{S0=0}[c][c]{\small $S_0 = 0\, $}
\psfrag{S1=1}[c][c]{\small $S_1 = 1$}
\psfrag{S2=2}[c][c]{\small $S_2 = 2$}
\psfrag{S3=4}[c][c]{\small $S_3 = 4$}
\psfrag{t0=0}[c][c]{\small $\tau_0 = 0$}
\psfrag{t1=4}[c][c]{\small $\tau_1 = 4$}
\psfrag{t2=7}[c][c]{\small $\tau_2 = 7$}
\psfrag{t3=8}[c][c]{\small $\tau_3 = 8$}
\psfrag{t4=11}[c][c]{\small $\tau_4 = 11$}
\psfrag{s0=0}[c][c]{\small $\sigma_0 = 0$}
\psfrag{s1=1}[c][c]{\small $\sigma_1 = 1$}
\psfrag{s2=5}[c][c]{\small $\sigma_2 = 5$}
\psfrag{s3=8}[c][c]{\small $\sigma_3 = 8$}
\psfrag{s4=9}[c][c]{\small $\sigma_4 = 9$}
\psfrag{Y1}[c][c]{\small $Y_1$}
\psfrag{Y2}[c][c]{\small $Y_2$}
\psfrag{Y3}[c][c]{\small $Y_3$}
\psfrag{Y4}[c][c]{\small $Y_4$}
\psfrag{Ds1}[c][c]{\small $\Delta_{\sigma_1}$}
\psfrag{Ds2}[c][c]{\small $\Delta_{\sigma_2}$}
\psfrag{Ds3}[c][c]{\small $\Delta_{\sigma_3}$}
\psfrag{Ds4}[c][c]{\small $\Delta_{\sigma_4}$}
\psfrag{Ds5}[c][c]{\small $\Delta_{\sigma_5}$}
\psfrag{Lb=}[c][c]{\small $L_b =$}
\includegraphics[width=14cm]{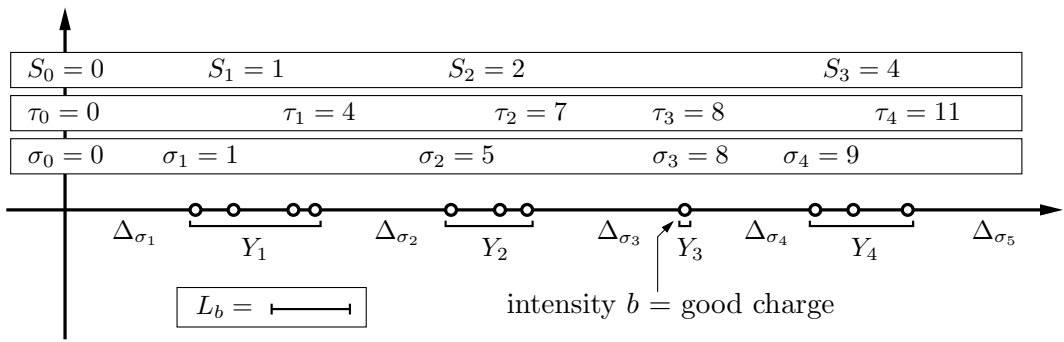}
\end{center}
\caption{A sample configuration of charges, with the corresponding
values of the variables $\gs_k$, $\tau_k$ and $S_k$.
Also indicated are the blocks $Y_k$ and the spacing $\Delta_{\gs_k}$
between blocks.}
\label{fig3}
\end{figure}

We define two sequences $\left\{  \sigma_{n}\right\}_{n\geq0}$, $\left\{
\tau_{n}\right\}  _{n\geq0}$ of random times by $\sigma_{0}:=0,$%
\[
\tau_{n}:=\inf\left\{  k\geq\sigma_{n}:\Delta_{k+1}>L_{b}\right\}\;,\qquad
\qquad \sigma_{n+1}=\tau_{n}+1 \,.
\]
Note that $\tau_{0}=0$ if and only if $\Delta_{1}>L_b=\big\lfloor \mathrm{e}^{\frac{2}%
{3}(b+K_{b})}\big\rfloor$. In words, the sequence $\{\tau_n\}_{n\ge 0}$
represents the indices of those positive charges that have a `distant'
($\Delta>L_b$) neighboring charge on the right.
Of course, one could define the sequence $\{\tau_n\}_{n\ge 0}$ alone,
without the need of introducing $\{\sigma_n\}_{n\ge 0}$, but in the sequel
it will be convenient to deal with both $\{\tau_n\}_{n\ge 0}$ and $\{\sigma_n\}_{n\ge 0}$.
Next we define a sequence of random variables $Y_{k},\ k\geq 0,$
taking values in the space%
\[
\Gamma:=\big(  \left\{  1\right\}  \times\mathbb{R}_{+}\big)  \cup
\bigcup\nolimits_{n\geq2} \big( \left\{  n\right\}  \times\mathbb{N}^{n-1}%
\times\mathbb{R}_{+}^{n} \big) \,,
\]
by
\[
Y_{k}:=\big( \, \tau_{k}-\sigma_{k}+1 \,,\, \left(  \Delta
_{\sigma_{k}+1},\ldots,\Delta_{\tau_{k}}\right) \,,\, \left(  \eta_{\sigma_{k}%
},\eta_{\sigma_{k}+1},\ldots,\eta_{\tau_{k}}\right)  \, \big)
\]
(the meaning of these definitions will be explained in a moment).
We occasionally write $\Gamma=\bigcup_{n=1}^{\infty}\Gamma_{n}.$ Here we
understand that in case $\tau_{k}=\sigma_{k},$ the $\Delta$-part is absent,
and the variable takes values in the $\Gamma_{1}$-part of $\Gamma.$ It should
be remarked that in case $\sigma_{0}=\tau_{0}=0$ we have $Y_{0}=\left(
1,0\right)$ (recall that $\eta_0 := 0$). See Figure~\ref{fig3} for a graphical
illustration of the variables just introduced.

Let us give some insight into these definitions.
Each $Y_k$ represents a block of adjacent positive charges, possibly reducing to
one single charge. More precisely:
\begin{itemize}
\item If $Y_k$ contains more than one charge, i.e.
$Y_k \in \Gamma_n$ with $n \ge 2$, then $Y_k$ corresponds
to a {\sl bad block} as defined in \S\ref{sec:renormalization},
because by construction the distances between the positive charges contained
in it, $\{\Delta_i\}_{\sigma_{k}+1\leq i\leq \tau_k}$, are such that
$\Delta_i\leq L_b$.
\item If on the other hand $Y_k \in \Gamma_1$,
then $Y_k$ is an {\sl isolated charge}, because
the distances $\Delta_{\sigma_{k}}, \Delta_{\sigma_{k+1}}$ with its
neighboring blocks are by construction larger than $L_{b}$.
\item If $Y_k \in \Gamma_1$ and in addition $\eta_{\gs_k} = b$,
i.e. $Y_k \in \{1\}\times\{b\}$, then $Y_k$ is a {\sl good charge}.
\end{itemize}
The variables $\sigma_k$ and $\tau_k$ give the indexes of the
first and last positive charge appearing in the $k$-th block, and therefore
$\tau_n-\sigma_n+1$ is the number of positive charges in a bad block.

\smallskip

Being larger than $L_b$, the variables $\Delta_{\sigma_{k}}$ are not geometrically distributed.
We therefore subtract $L_{b}$ and put%
\[
\widehat\Delta_{k}:=\Delta_{\sigma_{k}}-L_{b}\,, \qquad k\geq1\,.
\]
Note that the two sequences $\{Y_k\}_{k\ge 0}$ and $\{\widehat\Delta_k\}_{k\ge 1}$
contain all the information of the original sequence $\omega$.
The basic properties of the variables $Y_k$, $\widehat\Delta_k$ are
given in the following Lemma, whose (elementary) proof is omitted for conciseness.

\smallskip
\begin{lem}
\label{Le_Clump}$\,$
\begin{enumerate}
\item[a)] The random variables $Y_{0},Y_{1},Y_{2},\ldots,\widehat\Delta_{1}%
,\widehat\Delta_{2},\ldots,$ are independent.

\item[b)] The random variables $Y_{k},k\geq1,$ are identically distributed with the following
distribution:
\begin{itemize}
\item $\tau_{k}-\sigma_{k}+1$ is geometrically distributed with
parameter:%
\begin{equation}
q_{b}:={\mathbb{P}}_{\mu_{b}}(\Delta_{1}>L_{b})=(1-c_{b})^{L_{b}%
},\label{qb-Definition}%
\end{equation}
i.e. for $n \in \N$%
\[
P\left(  \tau_{k}-\sigma_{k}+1=n\right)  =q_{b}\, \left(  1-q_{b}\right)^{n-1} \,.
\]
\item Conditionally on $\{\tau_{k}-\sigma_{k}+1=n\}$, we have that $\eta_{\sigma_{k}%
},\eta_{\sigma_{k}+1},\ldots,\eta_{\tau_{k}}$ are $n$ i.i.d. random variables with
distribution $\overline{\mu}_{b}$ given by%
\[
\overline{\mu}_{b}\left(  x\right)  :=\frac{\mu_{b}\left(  x\right)  }{c_{b}%
}1_{\left\{  x\geq b\right\}  }.
\]
\item When $n\ge 2$, $\Delta_{\sigma_{k}+1},\ldots,\Delta_{\tau_{k}}$ are $n-1$ i.i.d.
random variables distributed like $\Delta_{1}$ conditionally on
$\left\{  \Delta_{1}\leq L_{b}\right\}$.
\end{itemize}

\item[c)] The random variables $\widehat \Delta_{1}, \widehat \Delta_{2},\ldots$ are i.i.d.
geometrically distributed with parameter $c_{b}.$

\item[d)] The distribution of $Y_{0}$ is given by an obvious modification:
$\tau_{0}+1$ is geometrically distributed as in b), and conditionally on
$\left\{  \tau_{0}+1=n\right\}  $ the distribution of the
$\Delta_{i},\eta_{i}$ is the same as described above, except that there is one
$\eta_{i}$ less, since $\eta_0 := 0$.
\end{enumerate}
\end{lem}
\smallskip

Next we define a mapping $\Phi:\Gamma\rightarrow \N \cup\{0\}.$ On $\Gamma_{1}$,
we simply put $\Phi\left(  \left(  1,\eta\right)  \right)  := \eta$, while on
$\Gamma_{n}$, $n\geq2$, we put%
\begin{equation}
\Phi\Big(  \big(  n,\left(  \Delta_{1},\ldots,\Delta_{n-1}\right)  ,\left(
\eta_{1},\ldots,\eta_{n}\right)  \big)  \Big)  \;:= \;\sum_{i=1}^{n}\eta
_{i}-\sum_{i=1}^{n-1}\left\lfloor \frac{3}{2}\log\Delta_{i}\right\rfloor
+2\left(  n-1\right)  K_{b}\,.\label{eq:size}%
\end{equation}
The interpretation is as follows: if $Y_k$ is a {\sl bad block}, i.e. if $Y_k \in \Gamma_n$
with $n\ge 2$, then $Y_k$ will be replaced by a single charge in the
new environment sequence $\omega'$, and $\Phi( Y_k)$ is
exactly the value of this {\sl clustered charge}. The reason why the
size of the clustered charge should be given by \eqref{eq:size}
will be clear in the next section.
Since $\eta_{i}\geq b$ and $\Delta_{i} \leq\mathrm{e}^{\frac{2}{3}(b+K_{b})}$,
it follows from \eqref{eq:size} that the value of the clustered charge
is always greater than $b + (n-1)K_b$, hence it is strictly greater
than $b$, if $n \ge 2$.

\begin{figure}[t]
\begin{center}
\psfragscanon
\psfrag{Y1}[c][c]{\small $Y_1$}
\psfrag{Y2}[c][c]{\small $Y_2$}
\psfrag{Y3}[c][c]{\small $Y_3$}
\psfrag{Y4}[c][c]{\small $Y_4$}
\psfrag{Ds1}[c][c]{\small $\Delta_{\sigma_1}$}
\psfrag{Ds2}[c][c]{\small $\Delta_{\sigma_2}$}
\psfrag{Ds3}[c][c]{\small $\Delta_{\sigma_3}$}
\psfrag{Ds4}[c][c]{\small $\Delta_{\sigma_4}$}
\psfrag{Ds5}[c][c]{\small $\Delta_{\sigma_5}$}
\psfrag{Lb=}[c][c]{\small $L_b =$}
\psfrag{go}[c][c]{$\omega$}
\psfrag{go'}[c][c]{$\omega'$}
\psfrag{D1}[c][c]{\small $\widehat\Delta_{1}$}
\psfrag{D2}[c][c]{\small $\widehat\Delta_{2}$}
\psfrag{D3+D4}[c][c]{\small $\widehat\Delta_{3} + \widehat\Delta_4$}
\psfrag{D5}[c][c]{\small \rule{0pt}{1.1em}$\qquad \quad \widehat\Delta_{5}$ (+\ldots)}
\psfrag{n1}[c][c]{\small $\eta'_{1}$}
\psfrag{n2}[c][c]{\small $\eta'_{2}$}
\psfrag{n3}[c][c]{\small $\eta'_{3}$}
\psfrag{Tb}[c][c]{$T_b$}
\includegraphics[width=14cm]{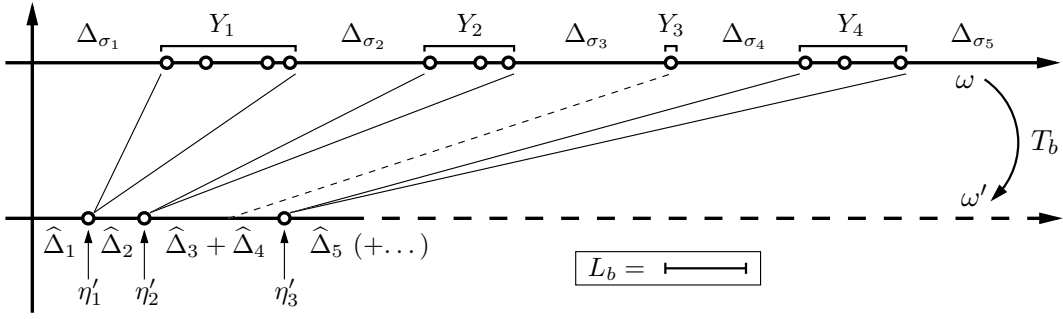}
\end{center}
\caption{An illustrative example of the map $T_b$ sending $\go$ to $\go'$
(the starting $\go$ is the same as in Figure~\ref{fig3}). The distances
$\Delta_{\gs_k}$ between blocks are shortened by $L_b$ and become $\widehat \Delta_k$,
the bad blocks ($Y_1, Y_2, Y_4$ in this example) are clustered into one single charge,
the good charges ($Y_3$ in this example) are erased, while the charges that
are isolated but not good (none in this example) are left unchanged.}
\label{fig4}
\end{figure}

We are ready to define the new sequence $\omega' = T_b(\omega)$.
First we set
\begin{equation} \label{eq:eta0}
    \widehat{\eta}_{0}:=\Phi\left(  Y_{0}\right)\,  .
\end{equation}
Then we introduce a sequence of stopping times by setting $S_0 := 0$ and for $k\ge 1$%
\[
S_{k} := \inf \big\{  n>S_{k-1}:\ \{ \tau_{n}>\sigma_{n} \} \ \mathrm{or\ }
\{\tau_n = \sigma_n \text{ and } \eta_{\sigma_{n}}>b\} \big\}
\]
(see also Figure~\ref{fig3}).
The variables $\{S_k\}_{k \ge 0}$ indicate which blocks will survive
after the renormalization: the block $Y_{S_k}$ will
become the $k$-th positive charge of $\go'$. More precisely,
we set
\begin{equation}\label{eq:omega'}
\eta_{k}^{\prime} \;:=\; \Phi\left(  Y_{S_{k}}\right)  \,, \qquad \qquad
\Delta_{k}^{\prime} \;:=\; \sum_{j=S_{k-1}+1}^{S_{k}}\widehat{\Delta}_{j} \,,
\end{equation}
and the sequence $\left(  \Delta_{k}^{\prime},\eta_{k}^{\prime}\right)_{k\geq1}$
defines our new sequence $\omega^{\prime}=:T_{b}\left(  \omega\right)$
(see also Figure~\ref{fig4}).
Note that the effect of the sequence $\{S_k\}_{k \ge 0}$ is to erase all the
{\sl good charges} in $\go$.

\smallskip

In the next lemma we show that $\{\go'_n\}_{n\ge 1}$ is indeed i.i.d.,
and we denote its one-marginal law by $\mu_{b+1} := T_b(\mu_b)$.
Observe that $\widehat{\eta}_{0}$ is not included in
the new sequence $\go'$, but it will enter the
estimate given below in Proposition \ref{prop:Zbound}
(remark that $\widehat{\eta}_{0}=0$ if $\tau_{0}=0$).

\smallskip
\begin{lem}
\label{Le_DistrRenorm}The random variables $\Delta_{k}^{\prime},\eta_{k}^{\prime}$
are all independent. The $\eta_{k}^{\prime}$ are identically distributed
and the $\Delta_{k}^{\prime}$ are identically geometrically distributed,
therefore the sequence $\{\omega'_n\}_{n\ge 1}$ is i.i.d.
The law $\mu_{b+1}(x)$ of $\omega'_1$ satisfies, for $x > 0$, the following relation:
\begin{equation}\label{eq:recursion}
\mu_{b+1}(x) \;=\; q_{b} \, \mu_{b}(x) \, \ind_{\{x\geq b+1\}}
\;+\; q_{b} \, \left(  Q_{b}\mu_{b}\right)  (x)\,,
\end{equation}
where%
\begin{equation}\label{eq:Qb}
\left(  Q_{b}\nu\right)  (x):=\sum_{n=1}^{\infty}\sum_{\ell_{1},\ldots
,\ell_{n}=1}^{L_{b}}\left(  \prod_{i=1}^{n}(1-c_{b})^{\ell_{i}-1}\right)
\sum_{\substack{x_{1},\ldots,x_{n+1}\geq b\\\sum_{i=1}^{n+1}x_{i}-\sum
_{i=1}^{n}\lfloor\frac{3}{2}\log\ell_{i}\rfloor+2nK_{b}=x}}\nu(x_{1})\cdots
\nu(x_{n+1})\,.
\end{equation}
\end{lem}
\smallskip

\proof
Introducing the subset $A$ of $\Gamma$ defined by
$A := \{1\} \times [b+1, \infty) \cup \bigcup_{n\ge 2} \Gamma_n$,
we see that $S_k$ is nothing but the sequence of times when $Y_n$
visits the set $A$. In particular, by Lemma~\ref{Le_Clump},
the increments $\{S_k - S_{k-1}\}_{k\ge 1}$ are i.i.d. geometrically
distributed, independent of the $\widehat\Delta_n$, with parameter $\gamma = P(Y_1 \in A)$.
(As a matter of fact, one can easily compute $\gamma = (1-q_b) + q_b \,{\tilde c}_{b}/c_b$,
cf. \eqref{cb-Definition},
but the precise value of $\gamma$ is immaterial for the proof.)
Therefore the variables $\Delta'_k$ are i.i.d. geometrically distributed,
with parameter $c_b \, \gamma$.

Notice that the variables $Y_{S_k}$ are i.i.d. and independent
of the $S_n$. It follows easily that the variables $\eta'_k$
are i.i.d. and independent from the $\Delta'_n$. The sequence $\omega'_n$ is
therefore i.i.d. and moreover for $x>0$ we have
$P(\omega'_1 = x) = (c_b \, \gamma) \cdot P(\eta'_1 = x)$.

It remains to determine the law of $\eta'_1 = \Phi(Y_{S_1})$.
Plainly, $Y_{S_1}$ is distributed like $Y_1$ conditionally on $\{Y_1 \in A\}$,
hence
\begin{equation} \label{eq:moreeas}
    P(\omega'_1 = x) \;=\; c_b \, \gamma \, P\big(\, \Phi(Y_1) = x \,\big|\, Y_1 \in A\,\big)
    \;=\; c_b \, P\big(\, \Phi(Y_1) = x \,, Y_1 \in A\,\big)\,.
\end{equation}
Since $\{Y_1 \in A\}$ is the disjoint union $\{\tau_1 = \gs_1,\, \eta_{\gs_1} > b\}
\cup \bigcup_{n\ge 1} \{\tau_1 - \gs_1 = n\}$ and since $\Phi(Y_1) = \eta_{\gs_1}$
if $\tau_1 = \gs_1$, we can write
\begin{equation*}
\begin{split}
    P(\omega'_1 = x) \;=\; & \; c_b \, P(\tau_1 = \gs_1) \, P\big(\, \eta_{\gs_1} = x,\,
    \eta_{\gs_1} > b \,\big|\, \tau_1 = \gs_1\,\big)\\
    & \;+\; \sum_{n=1}^\infty c_b \,
    P(\tau_1 - \gs_1 = n) \, P\big(\, \Phi(Y_1) = x \,\big|\, \tau_1 - \gs_1 = n\,\big)\,.
\end{split}
\end{equation*}
By Lemma~\ref{Le_Clump}, the first term in the r.h.s. equals
\begin{equation*}
    c_b \, q_b \, \ind_{\{x \ge b+1\}} \, \overline \mu_b(x) \;=\;
    q_b \, \ind_{\{x \ge b+1\}} \, \mu_b(x)\,,
\end{equation*}
matching with \eqref{eq:recursion}. Using \eqref{eq:size}
and again Lemma~\ref{Le_Clump}, we rewrite the second term as
\begin{equation*}
\begin{split}
    & \sum_{n=1}^\infty c_b \, q_b \, (1-q_b)^{n} \,
    \sum_{\ell_1, \ldots, \ell_n = 1}^{L_b} \, \sumtwo{x_1, \ldots, x_{n+1} \ge b}
    {\Phi(n+1,(\ell_i),(x_i)) = x} \, P\big( \, \widehat\Delta_{\gs_1 +i} = \ell_i,\,
    \eta_{\gs_1 + j - 1} = x_j,\,\big|\, \tau_1 - \gs_1 = n \big)\\
    & \quad \
    \;=\; \sum_{n=1}^\infty c_b \, q_b \, (1-q_b)^{n} \,
    \sum_{\ell_1, \ldots, \ell_n = 1}^{L_b} \, \sumtwo{x_1, \ldots, x_{n+1} \ge b}
    {\Phi(n+1,(\ell_i),(x_i)) = x} \, \prod_{i=1}^n \frac{c_b (1-c_b)^{\ell_i - 1}}{1-q_b}
    \,\prod_{j=1}^{n+1} \frac{\mu_b(x_j)}{c_b} \\
    & \quad \
    \;=\; q_b \,\sum_{n=1}^\infty \,
    \sum_{\ell_1, \ldots, \ell_n = 1}^{L_b} \, \sumtwo{x_1, \ldots, x_{n+1} \ge b}
    {\Phi(n+1,(\ell_i),(x_i)) = x} \, \prod_{i=1}^n (1-c_b)^{\ell_i - 1}
    \,\prod_{j=1}^{n+1} \mu_b(x_j) \;=\; q_b \, \left( Q_b \mu_b \right)(x)\,,
\end{split}
\end{equation*}
so that equation \eqref{eq:recursion} is proved.\qed
\medskip

As a side remark, we observe that by summing equation \eqref{eq:recursion}
(or, more easily, equation \eqref{eq:moreeas})
over $x \ge b+1 $ we obtain the following explicit formula for $\mu_{b+1}(0)$:
\begin{equation*}
    \mu_{b+1}(0) \;=\; 1 \;-\; \big(q_b \,\tilde c_{b} \,+\, (1-q_b) c_b \big)\,,
\end{equation*}
cf. \eqref{cb-Definition}.
Since $\tilde c_{b} \le c_b = 1 - \mu_b(0)$,
it follows that $\mu_{b+1}(0) \,\ge\, \mu_b(0)$,
i.e. at each renormalization step the density of positive charges decreases.


\medskip
\section{Renormalization of the partition function}
\label{sec:pf}

In the preceding section we have defined, at each level $b$
of the induction, a renormalizing map $T_b$ acting on the
environment sequence $\go$ and producing a renormalized sequence
$\go' = T_b(\go)$. In this section we show that, by replacing
$\go$ by $T_b(\go)$, one gets an upper bound on the free energy.
This will be the key to the proof of Theorem~\ref{th:main2} in
Section~\ref{sec:proof}. With some abuse of notation, we define
the map $T_b$ acting on the positive number $C$ by
\begin{equation}\label{eq:abuse}
    T_b(C) \;:=\; \big( 1 + B\,e^{-K_b}\, C \big) \cdot C\,,
\end{equation}
where $K_b$ is defined in \eqref{eq:Kb}
and $B$ is an absolute constant defined in Lemma~\ref{lem:ren_th} below.
Then we have the following

\smallskip
\begin{prop} \label{prop:Zbound}
There exists $b_{0}$ such that for every $b\geq b_{0}$,
for every $C \in (0, 2\cC]$ 
and for all $N\in\mathbb{N}$, there exists for ${\mathbb{P}}_{\mu
_{b}}$--a.e. $\omega$, a natural number $n\left(  \omega,N\right)  <\infty$
satisfying%
\begin{equation*}
N \;\leq\: n\left(  \omega,N\right)\,,  \qquad \quad t_{N}(T_{b}(\omega)) \;\leq\;
t_{n\left(  \omega,N\right)  }(\omega)\,,
\end{equation*}
and such that%
\begin{equation}
{\mathcal{Z}}_{n\left(  \omega,N\right)  }(\omega,C) \;\leq\; \mathrm{e}%
^{\hat{\eta}_{0}}\cdot{\mathcal{Z}}_{N}\big(  T_{b}(\omega),T_{b}%
(C)\big)  .\label{eq:Zbound}%
\end{equation}
\end{prop}
\smallskip

\proof
We set%
\[
n\left(  \omega,N\right)  :=\tau_{S_{N}}.
\]
The interpretation is as follows:
by construction (see Section~\ref{sec:env}) $Y_{S_N}$ is the $N$-th block
of charges of $\go$ that will survive after the renormalization,
and $n\left(  \omega,N\right)$ is the index
of the last positive charge in that block.
Therefore it is evident that $n(\go,N) \ge N$. Also
$t_{n\left(  \omega,N\right)  }(\omega) \ge t_{N}(T_{b}(\go))$ is
easy to check, because in the renormalization procedure leading from
$\go$ to $\go' = T_{b}(\go)$ the distances between charges are shortened
(see also Figure~\ref{fig4}).

For the rest of the proof, we fix $N\in\mathbb{N},$ and $n:=n\left(
\omega,N\right)  ,$ and we typically drop them from notations. The estimate is
purely deterministic and holds for any $\omega$ which has the property that
$n<\infty.$

We are going to work with the subsets of $\{0, t_1, t_2, \ldots, t_{n}\}$
and we need some notation.
We write $\mathcal{J}$ for the collection of intervals $I_{j}:=\left\{
t_{\sigma_{j}},\ldots,t_{\tau_{j}}\right\}  ,$ $1\leq j\leq S_{N}$.
Note that $\bigcup_{j=1}^{S_N} I_j = \{ t_{\gs_1}, \ldots, t_n \}$.
In $\cJ$ there are the $N$ `bad' intervals $\hat{I}%
_{j}:=I_{S_{j}},$ $1\leq j\leq N,$ and the other ones, we call
`good'. Note that the `good' intervals correspond to what we
have called {\sl good charges} in Sections~\ref{sec:strategy} and~\ref{sec:env},
while the `bad' intervals correspond to the {\sl bad blocks} and also
to the {\sl isolated charges} that are not good.
In particular, the good intervals are just single points: for this
reason, we also call them `good points'. The bad intervals may
be single points, too. We write
$\mathcal{J}^{\mathrm{bad}}$ for the set $\left\{  \hat{I}_{1},\ldots,\hat
{I}_{N}\right\}  $ of bad intervals, and $\mathcal{G}\subset\left\{
t_{1},\ldots,t_{n}\right\}  $ for the subset of good points.
The first interval $I_{0}=\left\{  0,\ldots,t_{\tau_{0}}\right\}$ is somewhat special. In case
$\tau_{0}=0,$ there is no charge (because $\omega_{0}=0$). In case $\tau
_{0}>0,$ this interval is of course `bad',
but we keep it separate from the others (remark that we don't take
it into $\mathcal{J}$).

We are in the situation where between the bad intervals, there may or may not be
good points in $\mathcal{G}.$ Also before the first bad interval,
i.e.  between $I_{0}$ and $\hat{I}_{1}$, there may or may not be
good points. If $X\subset\mathcal{J},$ we write $\mathcal{P}\left(
X\right)  $ for the set of subsets $F\subset\left\{  0,t_{1},\ldots
,t_{n}\right\}  $ which contain $0$ and $t_{n},$ and which have the property
that it has non void intersection with any interval in $X,$ and empty
intersection with any interval in $\mathcal{J}\backslash X.$ Then%
\begin{align*}
\mathcal{Z}_{n}\left(  \omega,C\right)    & =\sum_{X \subset \cJ,\, X\ni \hat I_N} \ \sum_{F\in
\mathcal{P}\left(  X\right)  }C^{\left\vert F\right\vert -1}\zeta\left(
F\right)  R\left(  F\right)  \\
& =\sum_{X \subset \cJ,\, X \ni \hat I_N} \ \sum_{F\in\mathcal{P}\left(  X\right)  }\zeta\left(  F\right)
\left(  C^{\left\vert F\cap I_{0}\right\vert -1}R\left(  F\cap I_{0}\right)
\right)  \prod\limits_{I\in X}\left(  C^{\left\vert F\cap I\right\vert
}R\left(  F\cap I\right)  \right)  ,
\end{align*}
where for a finite set $A=\left\{  s_{0},\ldots,s_{m}\right\}$, $m\ge 1$, we set%
\[
\zeta\left(  A\right)  :=\prod_{i=1}^{m}\left(  s_{i}-s_{i-1}\right)  ^{-3/2},
\]
and we put $\zeta(A) := 1$ in case $A$ reduces to a single point. We also set
\[
R\left(  F\right)  :=\exp\left[  \sum\nolimits_{x\in F}\omega_{x}\right]  ,
\]
where we recall that $\omega_{0}:=0.$ Note that the sum over $X$ is only over those which contain
the last interval in $\mathcal{J}$ (which is a bad one, namely $\hat I_N = I_{S_{N}}$),
in agreement with \eqref{eq:flast}.

The $\zeta\left(  F\right)  $ for $F\in\mathcal{P}\left(  X\right)  $ contains
the parts inside the $F\cap I,$ and the `interaction
part'. We want to split off this interaction part, and
estimate it by a bound which depends only on $X.$ If the intervals in $X$ are
(in increasing order) $I_{r_{1}},\ldots,I_{r_{k}}=I_{S_{N}},$ we write $M_{j}$
for the largest point of $I_{r_{j}},$ and $m_{j}$ for the smallest. Put%
\[
\hat{\zeta}\left(  X\right)  :=\left(  m_{1}-t_{\tau_{0}}\right)  ^{-3/2}%
\prod_{j=2}^{k}\left(  m_{j}-M_{j-1}\right)  ^{-3/2}.
\]
Then, if $F\in\mathcal{P}\left(  X\right)  $%
\[
\zeta\left(  F\right)  \leq\hat{\zeta}\left(  X\right)  \zeta\left(  F\cap
I_{0}\right)  \prod\limits_{I\in X}\zeta\left(  F\cap I\right)  .
\]
The inequality comes from the fact that if $x$ is the largest element of
$I\cap F,$ and $y$ is the smallest element of $I^{\prime}\cap F,$ $I$ an
interval below $I^{\prime},$ then $\left(  y-x\right)  ^{-3/2}\leq\left(
m^{\prime}-M\right)  ^{-3/2},$ $m^{\prime}$ being the smallest point in
$I^{\prime},$ and $M$ the largest of $I.$ We set for $I\in\mathcal{J}$%
\[
\alpha\left(  I\right)  :=\sum_{F\subset I,\ F\neq\emptyset}\left(
C^{\left\vert F\right\vert }R\left(  F\right)  \zeta\left(  F\right)  \right)
,
\]
and $\alpha\left(  I_{0}\right)  $ with one $C$-factor less. With this notation,
we have the estimate%
\[
\mathcal{Z}_{n}\left(  \omega,C\right)  \leq
\sum_{X \subset \cJ,\, X\ni \hat I_N} \hat{\zeta}\left(
X\right)  \alpha\left(  I_{0}\right)  \left[  \prod\limits_{I\in X}%
\alpha\left(  I\right)  \right]  .
\]

Next we claim that%
\begin{equation}
\alpha\left(  I\right)  \leq C\mathrm{e}^{2\left(  \left\vert I\right\vert
-1\right)  K_{b}}R\left(  I\right)  \zeta\left(  I\right)  =:C\exp\left[
\Phi_{I}\right]  ,\label{EstAlpha}%
\end{equation}
where $K_b$ is defined in \eqref{eq:Kb}. This
is evident if $I$ contains just one point, say $t_{j},$ in which case
$\alpha\left(  I\right)  =C\mathrm{e}^{\eta_{j}}.$ If $I$ contains more than
one point, then first observe that there are $2^{\left\vert I\right\vert
}-1\leq\mathrm{e}^{\left(  \left\vert I\right\vert -1\right)  K_{b}}$
possibilities to choose a non-empty subset $F\subset I$,
because $K_b \ge 2\log 2$ for $b \ge b_0$ and
$b_0$ large. Assume $I=\left\{  t_{\sigma},t_{\sigma
+1},\ldots,t_{\sigma+R}=t_{\tau}\right\}  ,$ and $F=\left\{  t_{j_{1}}%
,\ldots,t_{j_{m}}\right\}  ,$ so that%
\begin{equation}
C^{\left\vert F\right\vert }R\left(  F\right)  \zeta\left(  F\right)
=C^{m}\exp\left[  \sum\nolimits_{r=1}^{m}\eta_{j_{r}}\right]  \prod_{r=2}%
^{m}\left(  t_{j_{r}}-t_{j_{r-1}}\right)  ^{-3/2}.\label{EstAlpha2}%
\end{equation}
We can bound this from above by replacing $\left(  t_{j_{r}%
}-t_{j_{r-1}}\right)  ^{-3/2}$ by $\left(  t_{j_{r-1}+1}-t_{j_{r-1}}\right)
^{-3/2},$ and for the remaining gaps $t_{i+1}-t_{i}$, we simply use
$1\leq\mathrm{e}^{K_{b}}\mathrm{e}^{\eta_{i}}\left(  t_{i+1}-t_{i}\right)
^{-3/2}$ (we recall that $t_{i+1}-t_{i} \le e^{\frac 23(b + K_b)}$ 
and $\eta_i \ge b$ by construction).
Therefore, the right hand side of (\ref{EstAlpha2}) is bounded by%
\[
C^{m}\mathrm{e}^{K_{b}\left(  R-m+1\right)  }R\left(  I\right)  \zeta\left(
I\right)  \leq C\mathrm{e}^{\left(  \left\vert I\right\vert -1\right)  K_{b}%
}R\left(  I\right)  \zeta\left(  I\right),
\]
having used that $K_b \ge \log 2\cC \ge \log C$,
and (\ref{EstAlpha}) follows. Note that the definition of $\Phi_I$
matches with \eqref{eq:size}.

For the first interval $I_{0},$ we have a factor $C$ less on the right hand
side of (\ref{EstAlpha}), and therefore from (\ref{EstAlpha}) it follows that
$\alpha\left(  I_{0}\right)  \leq\mathrm{e}^{\hat{\eta}_{0}}$
(we recall that $\hat{\eta}_0$ is defined in \eqref{eq:eta0}). Combining, we
get%
\[
\mathcal{Z}_{n}\left(  \omega,C\right)  \leq\mathrm{e}^{\hat{\eta}_{0}}%
\sum_{X \subset \cJ,\, X\ni \hat I_N} \hat{\zeta}\left(  X\right)  C^{\left\vert X\right\vert }\left[
\prod\nolimits_{I\in X}\exp\left[  \Phi_{I}\right]  \right] \,.
\]
We have thus succeeded in clustering all the bad blocks.

\smallskip

Let us now fix $X^{\prime}\subset\mathcal{J}^{\mathrm{bad}}$,
$X'\ni \hat I_N$. Summing over $X$ with
$X\cap\mathcal{J}^{\mathrm{bad}}=X^{\prime}$ amounts to summing over all
subsets of $\mathcal{G}.$ Assume the intervals in $X^{\prime}$ are described
by the sequence%
\[
0\leq k_{0}<j_{1}\leq k_{1}<j_{2}\leq k_{2}<\cdots<j_{m}\leq k_{m}= n
\]
with $I_{0}=\left\{  t_{0},\ldots,t_{k_{0}}\right\}  ,$ $X^{\prime}=\left\{
\left\{  t_{j_{1}},\ldots,t_{k_{1}}\right\}  ,\ldots,\left\{  t_{j_{m}}%
,\ldots,t_{k_{m}}\right\}  \right\}$, $m=|X'|$. Write also $\mathcal{G}_{r},$ $1\leq
r\leq m$, for the set of good points between $t_{k_{r-1}}$ and $t_{j_{r}},$
and $\overline{\mathcal{G}}_{r}:=\mathcal{G}_{r}\cup\left\{  t_{k_{r-1}%
},t_{j_{r}}\right\}  .$ Then we can write the summation over $X$ with
$X\cap\mathcal{J}^{\mathrm{bad}}=X^{\prime}$ as the summation over
$A_{1}\subset\mathcal{G}_{1},\ldots,A_{m}\subset\mathcal{G}_{m}.$ For a single
$t_{j} \in\mathcal{G}$ we have $\Phi_{\left\{  t_{j}\right\}}=b.$ We therefore get%
\[
\sum_{X\in\cJ:\, X\cap\mathcal{J}^{\mathrm{bad}}=X^{\prime}}\hat{\zeta}\left(
X\right)  C^{\left\vert X\right\vert }\left[  \prod\nolimits_{I\in X}%
\exp\left[  \Phi_{I}\right]  \right]  =\left[  \prod\nolimits_{I\in X^{\prime
}}\exp\left[  \Phi_{I}\right]  \right]  \prod\limits_{r=1}^{m}\Xi\left(
b, C,\overline{\mathcal{G}}_{r}\right)  ,
\]
where for a finite subset $\mathcal{A=}\left\{  s_{0},s_{1},\ldots
,s_{m}\right\}  ,$ the $s_{i}$ ordered increasingly, we set%
\begin{equation}\label{eq:Xi}
\Xi\left(  b,C,\mathcal{A}\right)  \;:=\; \sumtwo{A \subset \cA}{A \supset\{s_0, s_m\}}
C^{|A|-1} \, e^{(|A|-2)b} \, \zeta(A)\,.
\end{equation}

Remark that the points in $\overline{\mathcal{G}}_{r}$ have inter-distances all
$> \mathrm{e}^{\frac23 \left(  b+K_{b}\right)} $.
We can therefore apply Lemma~\ref{lemRemoving} below, and obtain%
\[
\Xi\left(  b, C,\overline{\mathcal{G}}_{r}\right)  \leq\frac
{C \big( 1+BC\mathrm{e}^{-K_{b}} \big)}{\left(  t_{j_{r}}-t_{k_{r-1}}\right)  ^{3/2}}.
\]
We thus get%
\[
\mathcal{Z}_{n}\left(  \omega,C\right)  \;\leq\; \mathrm{e}^{\hat{\eta}_{0}}%
\sum_{X^{\prime}\subset\mathcal{J}^{\mathrm{bad}},\; X'\ni \hat I_N }\left[  \prod\nolimits_{I\in
X^{\prime}}\exp\left[  \Phi_{I}\right]  \right]  \prod\limits_{r=1}^{|X'|}%
\frac{C\left(  1+BC\mathrm{e}^{-K_{b}}\right)  }{\left(  t_{j_{r}}-t_{k_{r-1}%
}\right)  ^{3/2}}\,.
\]
Note that $C\left(  1+BC\mathrm{e}^{-K_{b}}\right)$ is by definition $T_b(C)$,
cf.~\eqref{eq:abuse}.

We are almost done. For the renormalized environment $\go' = T_b(\go)$,
defined in Section~\ref{sec:env}, we call $\{t'_k, \eta'_k\}_{k\in\N}$ the
locations and intensities of the positive charges of $\go'$ (see \S\ref{sec:notation}).
Consider the following correspondence: to
each bad interval $\hat I_l \in \cJ^{\mathrm{bad}}$ we associate the positive charge
$\eta'_l \in \go'$. Notice in fact that $\eta'_l = \Phi_{\hat I_l}$, see \eqref{EstAlpha}
and \eqref{eq:omega'}. Moreover, given two bad intervals $\hat I_l = \{ t_{j_l}, \ldots, t_{k_l} \}$,
$\hat I_m = \{ t_{j_m}, \ldots, t_{k_m} \} \in \cJ^{\mathrm{bad}}$, with $\hat I_l$ below $\hat I_m$,
we can bound $t_{j_m} - t_{k_l} \,>\, t'_m - t'_l$, because in passing from $\go$
to $\go'$ the distances $\Delta_{\gs_i}$ between intervals have been
shortened to $\widehat\Delta_i$ (it may be useful to look at Figure~\ref{fig4}).
Therefore we can bound $\mathcal{Z}_{n}\left(  \omega,C\right)$ by
\begin{equation*}
    \mathcal{Z}_{n}\left(  \omega,C\right) \;\le\;
    \mathrm{e}^{\hat{\eta}_{0}}
    \sumtwo{A =\{a_0, \ldots, a_k\} \subset \{0,\ldots, N\}}{0=a_0 < \ldots < a_k=N}\;
    \prod_{i=1}^{|A|} e^{\eta'_{a_i}} \, \frac{T_b(C)}{(t'_{a_i} - t'_{a_{i-1}})^{3/2}}
    \;=\; \mathrm{e}^{\hat{\eta}_{0}} \; \, \cZ_N \left(  \omega' , T_b(C) \right) \,,
\end{equation*}
where the last equality is just the definition \eqref{eq:start}
of the partition function, and the proof is completed.\qed

\medskip

We conclude this section with an auxiliary result (Lemma~\ref{lemRemoving} below)
that is used in the preceding proof. We first need a basic renewal theory lemma.

\smallskip
\begin{lem}
\label{lem:ren_th}
There exist positive constants $B$ and $K_{0}$ such that for
every $C>0$ and for all $K\geq K_{0}+\log C$ the following relation holds for
every $N\in\mathbb{N}$:
\[
\Theta_{N}^{+} \;:=\; \sum_{k=1}^{N}\sum_{\substack{j_{1},\ldots,j_{k-1}\in\mathbb{N}
\\0=:j_{0}<j_{1}<\ldots<j_{k-1}<j_{k}:=N}}\prod_{{\ell}=1}^{k}\frac
{C\mathrm{e}^{-K}}{(j_{{\ell}}-j_{{\ell}-1})^{3/2}}\leq\left(  1+BC\mathrm{e}%
^{-K}\right)  \,\frac{C\mathrm{e}^{-K}}{N^{3/2}}\,.
\]
\end{lem}
\smallskip

\proof
Defining the constant $A:=\left(  \sum_{n=1}^{\infty}n^{-3/2}\right)  ^{-1}$,
we set
\[
q(n):=\frac{A}{n^{3/2}}\,,\qquad \quad \gamma:=\frac{C}{A}\mathrm{e}^{-K}\,.
\]
Note that we can write
\begin{equation}
\Theta_{n}^{+} \;=\; \sum_{k=1}^{\infty}\gamma^{k}\,q^{\ast k}(n)\,,\label{eq:conv}%
\end{equation}
where $q^{\ast k}(\cdot)$ denotes the $k$-fold convolution of the probability
distribution $q(\cdot)$ with itself.

Let us prove by induction that $q^{\ast k}(n)\leq A\,k^{5/2}/n^{3/2}$ for
every $k\in\mathbb{N}$ and $n\in\mathbb{N}$. The case $k=1$ holds by
definition of $q(\cdot)$. For the inductive step, if $k$ is even, $k=2m$, we
can write
\begin{align*}
q^{\ast(2m)}(n) &  \leq2\sum_{\ell=1}^{\lfloor n/2\rfloor}q^{\ast m}%
(\ell)\,q^{\ast m}(n-\ell)\leq2\,A\,m^{5/2}\sum_{\ell=1}^{\lfloor n/2\rfloor
}q^{\ast m}(\ell)\,\frac{1}{(n-\ell)^{3/2}}\\
&  \leq\frac{2\,A\,m^{5/2}}{(n/2)^{3/2}}\sum_{\ell=1}^{\infty}q^{\ast m}%
(\ell)=\frac{A\,(2m)^{5/2}}{n^{3/2}}\,,
\end{align*}
and the odd case follows analogously. Then by \eqref{eq:conv} we can bound
$\Theta_{n}^{+}$ by $Ag\left(  \gamma\right)  n^{-3/2},$ where%
\[
g\left(  \gamma\right)  :=\sum_{k=1}^{\infty}k^{5/2}\gamma^{k}\leq
\gamma+8\gamma^{2},
\]
provided $0<\gamma\leq\gamma_{0},$ $\gamma_{0}$ sufficiently small. Let us set
$K_{0}:=-\log(A\,\gamma_{0})$, then if $K\geq K_{0}+\log C$, we have
$\gamma\leq\gamma_{0}$, and therefore
\[
\Theta_{n}^{+} \;\leq\; \left(  \gamma+8\,\gamma^{2}\right)  \frac{A}{n^{3/2}}\leq\left(
1+\frac{8C}{A}\mathrm{e}^{-K}\right)  \frac{C\mathrm{e}^{-K}}{n^{3/2}}.
\]
The proof is completed by setting $B:=8/A$.
\qed

\medskip

\smallskip
\begin{lem} \label{lemRemoving}
Let $K_{0}$ and $B$ be the constants of Lemma~\ref{lem:ren_th}.
Then $\forall N\geq2$, $\forall b>0$, $\forall C>0$, $\forall K\geq K_{0}+\log
C$ and for all $\cT \,=\,(t_{0},\ldots,t_{N})\,\in\,\mathbb{N}^{N+1}$ with
\[
t_{0}\,<\,t_{1}\,<\,\ldots\,<\,t_{N}\qquad\text{and}\qquad t_{n}%
-t_{n-1}\,>\,\mathrm{e}^{\frac{2}{3}(b+K)}\,,\ \forall n=1,\ldots,N\,,
\]
the following relation holds for $\Xi(b,C,\cT)$ (defined in \eqref{eq:Xi}):
\begin{align*}
\Xi(b,C,\cT) \;\leq \; \left(  1+BC\mathrm{e}^{-K}\right)  \frac{C}{\left(
t_{n}-t_{1}\right)  ^{3/2}}\,.
\end{align*}
\end{lem}
\smallskip

\proof
Expanding the definition of $\Xi(b,C,\cT)$ we can write:
\begin{equation*}
    \Xi(b,C,\cT) \;=\; e^{-b} \, \sum_{k=1}^{N}\sum_{\substack{j_{1}%
,\ldots,j_{k-1}\in\mathbb{N}\\0=:j_{0}<j_{1}<\ldots<j_{k-1}<j_{k}:=N}%
}\prod_{{\ell}=1}^{k}\frac{C\,\mathrm{e}^{b}}{(t_{j_{{\ell}}}-t_{j_{{\ell}-1}%
})^{3/2}}\,.
\end{equation*}
Fix a configuration $j_{1}<\cdots<j_{k-1}\subset\{1,\cdots,N-1\}$. Then we
have%
\begin{align*}
\prod_{{\ell}=1}^{k} &  \frac{1}{t_{j_{{\ell}}}-t_{j_{{\ell}-1}}}=\frac
{1}{t_{N}-t_{0}}\frac{\sum_{{\ell}=1}^{k}(t_{j_{{\ell}}}-t_{j_{{\ell}-1}}%
)}{\prod_{{\ell}=1}^{k}(t_{j_{{\ell}}}-t_{j_{{\ell}-1}})}\\
&  =\frac{1}{t_{N}-t_{0}}\sum_{{\ell}=1}^{k}\prod_{{\ell}^{\prime}%
\in\{1,\cdots,k\}\setminus{\ell}}\frac{1}{t_{j_{{\ell}^{\prime}}}-t_{j_{{\ell
}^{\prime}-1}}}\\
&  \leq\frac{\mathrm{e}^{-\frac{2}{3}(b+K)(k-1)}}{t_{N}-t_{0}}\sum_{{\ell}%
=1}^{k}\prod_{{\ell}^{\prime}\in\{1,\cdots,k\}\setminus{\ell}}\frac
{1}{j_{{\ell}^{\prime}}-j_{{\ell}^{\prime}-1}}\ \ \left[  \mathrm{since\ }%
t_{j_{\ell}}-t_{j_{\ell-1}}\geq\mathrm{e}^{2\left(  b+K\right)  /3}\left(
j_{\ell}-j_{\ell-1}\right)  \right] \\
&  =\frac{\mathrm{e}^{-\frac{2}{3}(b+K)(k-1)}}{t_{N}-t_{0}}\frac{\sum_{{\ell
}=1}^{k}(j_{{\ell}}-j_{{\ell}-1})}{\prod_{{\ell}=1}^{k}(j_{{\ell}}-j_{{\ell
}-1})}=\frac{\mathrm{e}^{-\frac{2}{3}(b+K)(k-1)}}{t_{N}-t_{0}}\frac{N}%
{\prod_{{\ell}=1}^{k}(j_{{\ell}}-j_{{\ell}-1})}.
\end{align*}
Therefore we get:
\begin{align*}
\Xi(b,C,\cT) & \; \leq\; e^{-b} \, \frac{\mathrm{e}^{b+K}\,N^{3/2}}%
{(t_{N}-t_{0})^{3/2}}\sum_{k=1}^{N}\sum_{\substack{j_{1},\ldots,j_{k-1}%
\in\mathbb{N} \\0=:j_{0}<j_{1}<\ldots<j_{k-1}<j_{k}:=N}}\prod_{{\ell}=1}%
^{k}\frac{C\,e^{-K}}{(j_{{\ell}}-j_{{\ell}-1})^{3/2}}\\
& \; \leq\; \left(  1+BC\mathrm{e}^{-K}\right)  \,\frac{C}{(t_{N}%
-t_{0})^{3/2}}\,,
\end{align*}
having used Lemma \ref{lem:ren_th}, and we are done.
\qed


\medskip
\section{Proof of Theorem \ref{th:main2}}
\label{sec:proof}

The starting point in the proof of Theorem~\ref{th:main2} is Proposition~\ref{prop:Zbound},
which immediately gives an upper bound on the free energy $\cf(\mu_b, C)$,
for every $b \ge b_0$ and $C \in (0, 2\cC]$. In fact,
since the limit in \eqref{eq:flast} holds $\bbP_{\mu_b}(\dd\go)$--a.s.,
we can take it along the (random) subsequence $\{n(\go,N)\}_{N\in\N}$ and
relation \eqref{eq:Zbound} yields
\begin{equation*}
    \cf \big(\mu_b, C \big) \;\le\; \liminf_{N \to \infty}
    \, \frac{1}{t_{n(\go, N)}(\go)} \, \log \cZ_N\big( T_b(\go), T_b(C) \big)\,,
    \qquad \quad \bbP_{\mu_b}(\dd\go)\text{--a.s.},
\end{equation*}
and since $t_{n(\go, N)}(\go) \ge t_{N}(T_b(\go))$ we obtain, again by \eqref{eq:flast},
\begin{equation*}
    \cf \big(\mu_b, C \big) \;\le\; \cf \big( \mu_{b+1}, T_b(C) \big)\,.
\end{equation*}
Note in fact that $\mu_{b+1}$ is by definition the one-marginal law of $T_b(\go)$,
when $\go$ has law $\mu_b$, see Section~\ref{sec:env}.
We now iterate this relation starting from $b=\gb$:
if we set
\begin{equation}\label{eq:Cb}
    \cC_\gb := \cC \qquad \text{and} \qquad
    \cC_b \;:=\; \big( T_{b-1} \cdot T_{b-2} 
    \cdots T_{\gb+1} \cdot T_\gb \big) (\cC) \quad
    \text{for $b > \gb$}\,,
\end{equation}
since in \S\ref{sec:Cb} below we show that $\cC_b \le 2\cC$ for every $b \ge \gb$,
provided $\gb$ is sufficiently large,
we can write
\begin{equation}\label{eq:basic_bound}
    \cf \big(\mu_\gb, \cC\big) \;\le\; \cf \big( \mu_b, \cC_b\big)\,,
    \qquad \quad \forall \gb \ge b_0,\ \forall b \ge \gb\,.
\end{equation}
We stress that, though not explicitly indicated,
both the law $\mu_b$ and the constant $\cC_b$ depend also on $\gb$,
which is the starting level of our procedure: however $\gb$ is kept
fixed in all our arguments. We recall that to prove Theorem~\ref{th:main2} it suffices
to show that $\cf \big(\mu_\gb, \cC\big) = 0$ when $\gb$ is large
(see \S\ref{sec:transformation}). Hence by \eqref{eq:basic_bound}
we are left with showing
that, if we fix $\gb$ sufficiently large, 
$\cf \big( \mu_b, \cC_b \big)$ vanishes as $b \to\infty$.

To estimate $\cf \big( \mu_b, \cC_b\big)$, we start from
a very rough upper bound on the partition function:
from the definition \eqref{eq:start} we can write for every $n\in\N$
\begin{equation*}
    \cZ_{n}(\go, \cC) \; \le \;
    e^{\sum_{i=1}^n (\eta_i + \log(\cC/A))}
    \left( \sum_{k=1}^{n} \sumtwo{j_1, \ldots, j_{k-1} \in \N}
    {0 =: j_0 < j_1 < \ldots < j_{k-1} < j_k = n} \prod_{\ell=1}^k
    \frac{A}{(t_{j_k} - t_{j_{k-1}})^{3/2}}\right) \,,
\end{equation*}
where $A := \big( \sum_{m=1}^\infty m^{-3/2} \big)^{-1} < 1$
is the constant that makes $m \mapsto A/m^{3/2}$ a probability law.
With this choice, the term in parenthesis in the r.h.s.
above is bounded from above by 
the probability that a renewal process with step law
$A/m^{3/2}$ visits the point $t_n$,
hence it is less than $1$ and we have
\begin{equation*}
    \cZ_{n}(\go, \cC) \; \le \;
    \exp \Bigg( \sum_{i=1}^n \bigg( \eta_i + \log\frac\cC A \bigg) \Bigg) \;=\;
    \exp \Bigg( \sum_{j=1}^{t_n(\go)} \bigg(\go_j + \log \frac \cC A \, \ind_{\{\go_j > 0\}} \bigg) \Bigg)
\end{equation*}
Now, if the sequence $\go$ is i.i.d. with marginal law $\mu$,
with $0 < \mu(0) < 1$, by \eqref{eq:flast} we have
\begin{align*}
    \cf(\mu, \cC) \;\le \; \lim_{n\to\infty} \, \frac{1}{t_n(\go)}
    \sum_{j=1}^{t_n(\go)}
    \left( \go_j + \log\frac{\cC}{A} \, \ind_{\{\go_j > 0\}} \right)
    \;=\; \bbE_\mu \big( \go_1 \big) + \bbP_\mu \big( \go_1 > 0 \big) \cdot \log\frac{\cC}{A}\,,
\end{align*}
having used that $t_n(\go) \to \infty$, $\bbP_\mu(\dd\go)$--a.s.,
and the strong law of the large numbers.
Combining this bound with \eqref{eq:basic_bound}, we get
\begin{equation}\label{eq:rough}
    \cf \big( \mu_\gb, \cC \big) \;\le\; \cf \big( \mu_b, \cC_b \big)
    \;\le\; \bbE_{\mu_b} (\go_1) \;+\; \big( 1 - \mu_b(0) \big) \cdot \log\frac{\cC_b}{A}\,,
\end{equation}
for every $\gb \ge b_0$ and for every $b \ge \gb$.

We are left with estimating the r.h.s. of \eqref{eq:rough}. To this
purpose, we exploit the stochastic domination on $\mu_b$ given by the following

\smallskip
\begin{prop}\label{propInduction}
There exists a finite $b_1$ (depending on $\cC > 0$ and $c > \frac 23$) such that for all
$b\ge \gb \ge b_1$ we have:
\begin{equation} \label{eq:stocdom}
\mu_b(x)\leq e^{-\frac 23 x-\sqrt{x}} \qquad \quad \forall x \ge b\,.
\end{equation}
\end{prop}
\smallskip

\noindent
Applying \eqref{eq:stocdom} to \eqref{eq:rough}, for any fixed $\gb \ge \gb_0 := b_0 \vee b_1$,
we obtain for every $b \ge \gb$
\begin{equation*}
    \cf \big( \mu_\gb, \cC \big) \;\le\; \sum_{x=b}^\infty x \, e^{-\frac 23 x - \sqrt x}
    \;+\; \left(\log\frac{\cC_b}{A}\right) \, \sum_{x=b}^\infty e^{-\frac 23 x - \sqrt x}\,.
\end{equation*}
It is clear that both the sums in the r.h.s. can be made arbitrarily small by
taking $b$ large. Moreover, in \S\ref{sec:Cb} below we show that $\cC_b \le 2 \cC$
for all $b$. Therefore, by letting $b\to\infty$,
we have shown that $\cf \big( \mu_\gb, \cC \big) = 0$ for all $\gb \ge \gb_0$,
and this completes the proof of Theorem~\ref{th:main2}.\qed

\bigskip

The proof of Proposition \ref{propInduction} is given in \S\ref{sec:induction}
below. Before that, we need to establish two technical lemmas.

\smallskip
\begin{lem}\label{lemInduction1}
Let $m\in\NN,\,m\geq 2$, let $b\in\N,\,b\geq \gb \vee 100$, and
$z\geq mb.$ Then%
\[
A_{m,b}\left( z\right) \;:=\; \sum_{\substack{ %
x_{1},\ldots ,x_{m}\geq b \\ x_{1}+\cdots +x_{m}=z}}\exp \left[
-\sum\nolimits_{i=1}^{m}\sqrt{x_{i}}\right] \leq \mathrm{e}^{-\sqrt{z}-\frac{%
m-1}{4}\sqrt{b}}.
\]
\end{lem}
\smallskip

\proof
We first treat the case $m=2.$ There we have for the left hand side%
\[
A_{2,b}\left( z\right) \le
2\sum_{x=b}^{\lceil z/2 \rceil}\exp \left[ -\sqrt{x}-\sqrt{z-x}%
\right] .
\]%
For $x\leq z/2$%
\[
\sqrt{z-x}=\sqrt{z}\sqrt{1-\frac{x}{z}}.
\]%
For $t\leq 1/2,$ we have from concavity of the square root function%
\[
\sqrt{1-t}\geq 1-\left( 2-\sqrt{2}\right) t,
\]%
and therefore for $x\leq z/2$%
\begin{eqnarray*}
\sqrt{z-x} &\geq &\sqrt{z}-\left( 2-\sqrt{2}\right) \frac{x}{\sqrt{z}} \\
&\geq &\sqrt{z}-\left( \sqrt{2}-1\right) \sqrt{x}.
\end{eqnarray*}%
Hence,%
\[
\sum_{x=b}^{z-b}\exp \left[ -\sqrt{x}-\sqrt{z-x}\right] \leq 2 \mathrm{e}^{-%
\sqrt{z}}\sum_{x=b}^{\infty }\exp \left[ -\frac{1}{2}\sqrt{x}\right] .
\]%
\begin{eqnarray*}
\sum_{x=b}^{\infty }\exp \left[ -\frac{1}{2}\sqrt{x}\right]  &\leq
&\int_{b-1}^{\infty }\mathrm{e}^{-\sqrt{t}/2}dt=8\int_{\sqrt{b-1}/2}^{\infty
}t\mathrm{e}^{-t}dt \\
&=&8\left( \frac{\sqrt{b-1}}{2}+1\right) \mathrm{e}^{-\sqrt{b-1}/2}\leq \frac 12 \exp %
\left[ -\frac{1}{4}\sqrt{b}\right] ,
\end{eqnarray*}%
if $b\geq 100$. This proves the claim for $m=2.$

The general case follows by induction on $m.$ Assume $m\geq 3.$ Then%
\begin{eqnarray*}
A_{m,b}\left( z\right)  &=&\sum_{x_{1}=b}^{z-\left( m-1\right) b}\mathrm{e}%
^{-\sqrt{x_{1}}}\sum_{\substack{ x_{2},\ldots ,x_{m}\geq b \\ x_{2}+\cdots
+x_{m}=z-x_{1}}}\exp \left[ -\sum\nolimits_{k=2}^{m}\sqrt{x_{k}}\right]  \\
&\leq &\sum_{x_{1}=b}^{z-\left( m-1\right) b}\mathrm{e}^{-\sqrt{x_{1}}}\exp %
\left[ -\sqrt{z-x_{1}}-\frac{m-2}{4}\sqrt{b}\right]  \\
&\leq &\mathrm{e}^{-\frac{m-2}{4}\sqrt{b}}\sum_{x=b}^{z-b}\exp \left[ -\sqrt{%
x}-\sqrt{z-x}\right] \leq \mathrm{e}^{-\sqrt{z}-\frac{m-1}{4}\sqrt{b}},
\end{eqnarray*}%
by induction and the $m=2$ case.
\qed

\medskip

\smallskip
\begin{lem}\label{lemInduction2}
There exists $b_2=b_2(\cC)$ such that for every $b\ge b_2$ and $n\in\N$ we have:
\begin{align*}
B_{n,b}(x) &\;=\; \sum_{\l_1,\cdots,\l_n=1}^{\lfloor e^{\frac{2}{3}(b+K_b)} \rfloor}
e^{-\frac{2}{3}(x+\sum_{i=1}^n
(\lfloor \frac{3}{2}\log\l_i \rfloor -2K_b))} A_{n+1,b}\left(x+\sum_{i=1}^n
\left(\left\lfloor \frac{3}{2}\log\l_i \right\rfloor
-2K_b\right)\right)\\
& \;\leq\; e^{-\frac{2}{3}x-\sqrt{x}}e^{-\frac{n}{8}\sqrt{b}},
\end{align*}
where $A_{m,b}(x)$ is defined in Lemma~\ref{lemInduction1} and $K_b$
is defined by \eqref{eq:Kb}.
\end{lem}
\smallskip

\proof
By Lemma \ref{lemInduction1}, we have:
\begin{equation*}
 A_{n+1,b}\left(x+\sum_{i=1}^n
\left( \left\lfloor\frac{3}{2}\log\l_i \right\rfloor
-2K_b\right)\right)\leq e^{-\frac{n}{4}\sqrt{b}}
e^{-\sqrt{x+\sum_{i=1}^n(\lfloor \frac{3}{2}\log\l_i \rfloor -2K_b)}}.
\end{equation*}
Since for $0 \le a \le b$ we have $\sqrt{b-a} \ge \sqrt b - \sqrt a$,
it follows that
\begin{equation*}
\sqrt{x+\sum_{i=1}^n \left( \left\lfloor \frac{3}{2}\log\l_i \right\rfloor-2K_b \right)} \ge
\sqrt x - \sqrt{\sum_{i=1}^n \left(2K_b - \left\lfloor \frac{3}{2}\log\l_i \right\rfloor \right)}
\ge \sqrt{x} -\sqrt{2nK_b}\,,
\end{equation*}
because $\ell_i  \ge 1$. Therefore
\begin{equation*}
 A_{n+1,b}\left(x+\sum_{i=1}^n
\left(\left\lfloor \frac{3}{2}\log\l_i \right\rfloor -2K_b\right)\right)
\leq e^{-\frac{n}{4}\sqrt{b}}e^{-\sqrt{x}}e^{\sqrt{2nK_b}}.
\end{equation*}
Plugging this in the definition of $B_{n,b}(x)$
and using the fact that $\lfloor \frac 32 \log \ell_i \rfloor \ge \frac 32 \log \ell_i -1$ yields:
\begin{align*}
B_{n,b}(x)&\leq
e^{-\frac{2}{3}x-\sqrt{x}}e^{-\frac{n}{4}\sqrt{b}}e^{\frac{4}{3}nK_b+\sqrt{2nK_b}}
\left(\sum_{\l_1=1}^{\lfloor e^{\frac{2}{3}(b+K_b)} \rfloor}\frac{e}{\l_1}
\right)^n\\
&\leq e^{-\frac{2}{3}x-\sqrt{x}}e^{-\frac{n}{4}\sqrt{b}}
e^{\frac{4}{3}n K_b+\sqrt{2nK_b}}\left(e \cdot \frac{2}{3}(b+K_b)\right)^n\\
&=e^{-\frac{2}{3}x-\sqrt{x}}e^{-\frac{n}{8}\sqrt{b}}
e^{(-\frac{1}{8}\sqrt{b}+\frac{4}{3}K_b+
\frac{\sqrt{2K_b}}{\sqrt{n}}+ 1 +\log \frac{2}{3}(b+K_b))n}.
\end{align*}
Observe that, for every fixed $\cC > 0$, we have
$-\frac{1}{8}\sqrt{b}+\frac{4}{3}K_b+\frac{\sqrt{K_b}}{\sqrt{n}}+ 1 + \log
\frac{2}{3}(b+K_b) \to -\infty$ as $b\rightarrow\infty$ (recall the
definition \eqref{eq:Kb} of $K_b$).
Therefore there exists $b_2 = b_2(\cC)$, such that for all $b\geq b_2$,
$B_{n,b}(x)\leq e^{-\frac{2}{3}x-\sqrt{x}}e^{-\frac{n}{8}\sqrt{b}}$.\qed


\subsection{Proof of Proposition~\ref{propInduction}}
\label{sec:induction}

The law $\mu_b$ is defined recursively by equations
\eqref{eq:recursion} and \eqref{eq:Qb},
starting from $\mu_\gb = (1-e^{-c\,\gb}) \gd_{\{0\}} + e^{-c\,\gb} \gd_{\{\gb\}}$.
In particular we have
\begin{equation*}
    \mu_\beta(\beta) \;=\; e^{-\beta c},\qquad \mu_\beta(x) \;=\; 0,\quad \forall\,x\geq\beta+1\,,
\end{equation*}
and for $x \ge \gb$
\begin{align} \label{eq:inductive}
\mu_{b+1}(x) \;\le\; \mu_b(x)+
\sum_{n=1}^\infty\;
\sum_{\l_1,\cdots,\l_n=1}^{\lfloor e^{\frac{2}{3}(b+K_b)} \rfloor}
\sumtwo{x_1,\cdots,x_{n+1}\geq b}{
\sum_{i=1}^{n+1}x_i-\sum_{j=1}^n(\lfloor \frac{3}{2}\log \l_j \rfloor - 2K_b ) = x}
\mu_{b}(x_1)\cdots\mu_b(x_{n+1}) \,.
\end{align}
We recall that $\mu_b$ is supported on $\{0\} \cup \{b,b+1,\ldots\}$.

Let us prove equation \eqref{eq:stocdom} by induction on
$b=\beta,\beta+1,\cdots$. For $b=\gb$ we have
\begin{equation*}
\mu_{\beta}(\beta)=e^{-\beta c}=e^{-\frac{2}{3}\beta-(c-\frac 23)\beta}\leq
e^{-\frac{2}{3}\beta-\sqrt{\beta}},
\end{equation*}
whenever $\beta \geq b_3(c) := \frac{1}{(c-\frac 23)^2}$. In particular,
\eqref{eq:stocdom} holds for $b=\gb$, provided $\beta\geq b_3(c)$.

Assume now that $\mu_\alpha(x)\leq e^{-\frac 23 x-\sqrt{x}}$, for all
$\beta\leq\alpha\leq b$. Then, using \eqref{eq:inductive}
and the fact that $\mu_\beta(x)=0$, $\forall x\geq
\beta+1$, we have:
\begin{equation*}
\mu_{b+1}(x) \;\le\; \sum_{\alpha=\beta}^b \sum_{n=1}^\infty
\sum_{\l_1,\cdots,\l_n=1}^{\lfloor e^{\frac{2}{3}(\alpha+K_b) \rfloor}}
\sum_{\scriptsize
\begin{array}{c}
x_1,\cdots,x_{n+1}\geq \alpha\\
\sum_{i=1}^{n+1}x_i-\sum_{j=1}^n(\lfloor \frac{3}{2}\log \l_j \rfloor -2K_b)=x
\end{array}
}
\mu_{\alpha}(x_1)\cdots\mu_\alpha(x_{n+1})\,,
\end{equation*}
Plugging in the induction assumption yields:
\begin{equation*}
\mu_{b+1}(x)\leq \sum_{\alpha=\beta}^b \sum_{n=1}^\infty B_{n,\alpha}(x),
\end{equation*}
where $B_{n,\alpha}(x)$ is defined in Lemma~\ref{lemInduction2}.
Assuming $\beta\geq b_2$ and using Lemma \ref{lemInduction2} gives:
\begin{align*}
\mu_{b+1}(x)&\leq e^{-\frac{2}{3}x-\sqrt{x}}\sum_{\alpha=\beta}^b
\sum_{n=1}^\infty e^{-\frac{n}{8}\sqrt{\ga}}\\
&\leq   e^{-\frac{2}{3}x-\sqrt{x}}\sum_{\alpha=\beta}^\infty
e^{-\frac{\sqrt{\ga}}{10}}\\
&\leq e^{-\frac{2}{3}x-\sqrt{x}}\int_{\beta-1}^\infty
e^{-\frac{\sqrt{\ga}}{10}}
= e^{-\frac{2}{3}x-\sqrt{x}}
\left(20(10+\sqrt{\beta-1})e^{-\frac{\sqrt{\beta-1}}{10}}\right)\\
&\leq e^{-\frac{2}{3}x-\sqrt{x}}, \mbox{ when $\beta\geq 6000$}.
\end{align*}
Hence Proposition \ref{propInduction} holds for all
$\beta\ge b_1(\cC,c) := \max \{b_2(\cC), b_3(c) ,6000\}$.
\qed
\medskip


\subsection{Bounding the constant $\cC_b$}
\label{sec:Cb}

We are going to show that $\cC_b \le 2 \cC$ for all $b \ge \gb$. We recall
the definition \eqref{eq:iter_bound} of $\cC_b$ (see also \eqref{eq:Cb}):
\begin{equation*}
    \cC_\gb \;:=\; \cC\,, \qquad \qquad
    \cC_b \;=\; T^{b-1} (\cC_{b-1}) \;=\;
    \big( 1 + B \, \cC_{b-1} \, e^{-K_{b-1}} \big) \cdot \cC_{b-1}\,,\qquad
    b > \gb\,,
\end{equation*}
where $\cC$ is the constant appearing in the definition
\eqref{eq:start} of the partition function.
Passing to logarithms we get
\begin{equation*}
    \log \cC_b \;=\; \log \cC  +
    \sum_{a = \gb}^{b-1} \log \big( 1 + B \, \cC_{a} \, e^{-K_a} \big)\,.
\end{equation*}

We prove that $\cC_b \le 2\cC$ by induction: the case $b=\gb$ is trivial since $\cC_\gb = \cC$.
Assuming that $\cC_a \le 2\cC$ for all $a \in \{\gb, \ldots, b-1\}$, we have
\begin{align*}
    \log \, \cC_b \;\le\; \log\cC \;+\; \sum_{a=\gb}^{b-1}
    \log \big( 1 + 2\, B \, \cC \, e^{-K_a} \big)\,.
\end{align*}
By the definition \eqref{eq:Kb} of $K_a$
we have $K_a \ge \big( K_0 + \log(2\cC) + 2 \log a\big) - 1$, therefore
\begin{align*}
    \log \, \cC_b \;\le\; \log\,\cC \;+\; \sum_{a=\gb}^{b-1}
    \log \big( 1 + e^{1-K_0} \frac{B}{a^2} \big)
    \;\le\; \log\cC \;+\; B\,e^{1-K_0} \sum_{a=\gb}^\infty \frac{1}{a^2}\,.
\end{align*}
Therefore, if we choose $\gb$ sufficiently large, we get $\log \cC_b \,\le\, \log(2\cC)$
and we are done.


\medskip
\section*{Acknowledgments}

We are very grateful to Giambattista Giacomin and
Ofer Zeitouni, for sharing with us their experience of the problem,
and to Dimitris Cheliotis for valuable observations on the manuscript.
This work is partially supported by the Swiss National Science Fundation
under contracts 200020--116348 (E.B.), 201121--107890/1 (F.C.) and 47102009 (B.d.T.).



\end{document}